\documentclass[11pt,e-only]{amsart}
\usepackage[english]{babel}
\usepackage{amsmath, amsfonts, amssymb, amsthm}
\usepackage{mathtools}
\usepackage{graphicx}
\usepackage{xcolor}
\usepackage{makecell}
\usepackage{stackengine}
\usepackage{longtable}
\renewcommand{\arraystretch}{1.2}
\usepackage{aliascnt}
\usepackage[separate-uncertainty = true, detect-all]{siunitx}
\usepackage{booktabs}
\usepackage{tikz}
\usetikzlibrary{braids,decorations.pathreplacing}
\usetikzlibrary{matrix,arrows.meta}
\usepackage{multirow}
\usepackage{array}
\usepackage[breaklinks=true]{hyperref}
\usepackage{xurl}
\hypersetup{
    colorlinks,
    linkcolor={blue!80!black},
    citecolor={blue!50!black},
    urlcolor={blue!80!black}
}
\usepackage[capitalize]{cleveref}
\usepackage[margin=3cm]{geometry}
\usepackage{caption}
\usepackage{setspace}
\usepackage{fancyhdr}
\usepackage{subcaption}
\usepackage[normalem]{ulem}
\usepackage[isbn=false,giveninits=true,date=year]{biblatex}
\usepackage{csquotes}

\renewbibmacro*{journal}{%
  \iffieldundef{shortjournal}
    {%
      \iffieldundef{journaltitle}
        {}
        {%
          \printtext[journaltitle]
            {%
              \printfield[titlecase]{journaltitle}%
              \setunit{\subtitlepunct}%
              \printfield[titlecase]{journalsubtitle}%
             }%
         }%
    }
    {\printtext[journaltitle]{\printfield[titlecase]{shortjournal}}}%
}

\renewbibmacro{in:}{}
\DeclareFieldFormat{url}{%
  \href{https://#1}{#1}%
}

\DeclareFieldFormat{doi}{%
  \href{https://doi.org/#1}{#1}%
}

\DeclareFieldFormat{eprint:arxiv}{%
  \href{https://arxiv.org/abs/#1}{#1}%
}

\AtEveryBibitem{%
  \iffieldundef{doi}{}{\clearfield{url}}%
}
\AtEveryBibitem{%
  \iffieldundef{eprint}{}{\clearfield{url}\clearfield{urldate}}%
}

\AtBeginBibliography{\setlength{\emergencystretch}{2em}}

\usepackage{mymacros}

\stackMath

\title{Quantum invariants and fiberedness}

\author[P. Orland]{Paul Orland}

\author[L. San Mart\'in Su\'arez]{Lara San Mart\'in Su\'arez}

\author[T. Saunders-A'Court]{Toby Saunders-A'Court}

\author[J. Svoboda]{Josef Svoboda}

\thanks{California Institute of Technology, Pasadena, CA 91125, USA}
\thanks{Emails: \texttt{porland@caltech.edu,lsanmart@caltech.edu,tsaunder@caltech.edu,svo@caltech.edu}}

\date{\today}

\begin{document}

\begin{abstract}
We explore the topological significance of the Gukov--Manolescu knot series $F_K$. We show that the leading coefficient of $F_K$ is a monomial and express its exponent in terms of the Hopf invariant for all homogeneous braid knots, and for fibered knots up to 12 crossings. As an application, we deduce an explicit formula for the Hopf invariant in terms of colored Jones polynomials. For non-fibered strongly quasipositive knots, we study a relation between $F_K$ and the stability series of the colored Jones function, and explore similarities between $F_K$ and knot Floer homology. Finally, we propose a slope conjecture for $F_K$, relating it to the boundary slopes of the knot.
\end{abstract}

\maketitle
\section{Introduction}

What do quantum invariants say about the topology and geometry of manifolds? To name just a few examples, the Jones polynomial \cite{Jones} helped resolve several questions in knot theory, including the famous Tait conjectures on alternating knots \cite{Murasugi1987, Kauffman1987, Thistlethwaite1987}. Its colored version \cite{RT} is related to hyperbolic invariants of knots and 3-manifolds, and to boundary slopes of essential surfaces \cite{garoufalidis2010slopes}. Khovanov homology \cite{Kh}, a categorification of the Jones polynomial, has led to numerous topological applications, including Piccirillo's proof that the Conway knot is not slice \cite{Piccirillo2020}.

The colored Jones polynomials have a 3-manifold analog: the Witten--Reshetikhin--Turaev (WRT) invariants, defined at roots of unity \cite{Witten1989, RT}. Ever since Khovanov's categorification of the Jones polynomial, the quest for a theory categorifying the WRT invariants has been a major open problem in quantum topology. In order to tackle this question, multiple new quantum invariants have appeared recently; among them the Gukov--Putrov--Vafa--Pei (GPPV) invariants \cite{GPPV}, which are analytic functions on the unit disc, recovering the WRT invariants as radial limits \cite{Murakami2024}. GPPV invariants were originally defined for closed negative definite plumbed 3-manifolds. In a search for a more general definition, Gukov and Manolescu in \cite{GM} proposed a new knot invariant $F_K$, as a knot-theoretic counterpart of GPPV invariants, and defined $F_K$ for certain plumbed knots. In \cite{Park20,Park21}, Park extended the $F_K$ invariant to certain classes of hyperbolic knots using quantum groups and Verma modules. The key idea is his \emph{inverted state sum}, which uses an extension of the $R$-matrix, along with certain additional data on a braid representative of the knot.

The goal of this paper is to demonstrate that $F_K$ detects interesting topological quantities such as fiberedness, the Hopf invariant of fibered knots, and boundary slopes of essential surfaces, and explore the relations between $F_K$ and other knot invariants. Our findings suggest that $F_K$ behaves quite distinctly from the colored Jones polynomial, and in fact, inherits some characteristics similar to knot Floer homology. This suggests that the (hypothetical) categorification of $F_K$ and GPPV invariants could have interesting properties and lead to new topological applications.

\subsection{GM series and colored Jones function}

In this paper, we consider the quantum invariants associated with the quantum group $U_q(\sl_2)$ (of root system $A_1$). The \emph{colored Jones function} $J_K$ of a knot $K$ is a sequence of Laurent polynomials $J^{K}_{n}$ for $n=1,2,\dots$ They are defined as the Reshetikhin--Turaev invariants associated to the $n$-dimensional representation of $U_q(\sl_2)$. Recall two fundamental properties of the colored Jones function $J_K$:
\begin{itemize}
    \item $J_K$ is $q$-holonomic, i.e., it satisfies a $q$-recurrence relation defined by the \emph{noncommutative $A$-polynomial} $\widehat{A}$ of the knot $K$ \cite{GaroufalidisLe05}, which is conjecturally given by the quantization of the classical $A$-polynomial; cf. \cite{Garoufalidis2004, Guk05}. More specifically, the sequence $J^K_n$ is annihilated by an element $\widehat{A}$ of the Weyl algebra $\Z[q^{\pm 1}]\langle Q,E \rangle/EQ - qQE$, where the operator $Q$ multiplies $J^K_n$ by $q^n$ and the operator $E$ shifts the index $n$ by 1.

    \item $J_K$ admits the Melvin--Morton--Rozansky (MMR) expansion near $q=1$, whose coefficients are rational functions, with denominators given by powers of the Alexander polynomial of the knot $K$ \cite{MM,BNG,RozMMRConj}. 
\end{itemize}

The GM series, or $F_K$ invariant, of a knot $K$ in $S^3$ is a formal series of the form
\begin{equation}\label{eq:F_K-terms}
    F_K(x,q) = x^{\frac12} \sum_{n=0}^{\infty} f_n(q) x^{d_K-1+n},
\end{equation}
where $f_n(q)$ are Laurent power series with integer coefficients, with $f_0(q) \neq 0$ and $d_K \in \Zpos$.

The existence of the $F_K$ invariant is still conjectural in general. However, it was defined by Park for a large class of knots, called \emph{nice} (see \cref{def:nice-knot}), which includes homogeneous braid knots. For nice knots, Park proved that $F_K$ satisfies the above properties of the colored Jones function $J_K$, with important differences:

\begin{itemize}
    \item The coefficient sequence $(f_n(q))_{n \geq 0}$ satisfies a $q$-recurrence given by the `transpose' of the noncommutative A-polynomial, i.e., by $\widehat{A}$ with $E$ replaced by $Q$ and vice versa.
    \item The $F_K$ admits an expansion at $q=1$ whose coefficients are rational functions in $x$. This series coincides with the MMR expansion of $J_K$, but rational functions are expanded at $x=0$, rather than at $x=1$.
\end{itemize}

This is as far as the similarities between $J_K$ and $F_K$ go. Despite the fact that the $F_K$ invariant was largely motivated by $J_K$ and WRT invariants (and GPPV invariants), it seems that the properties of the $F_K$ invariant are quite different and we are at the beginning of a completely new and rich theory.

\subsection{Main results} The primary focus is the first nonzero coefficient $f_0(q)$. Our first result states that the minimal coefficient $f_0$ of a nice knot $K$ is a monomial. More significantly, for a large class of fibered knots---including all homogeneous braid knots \cite{Stallings1978} and fibered knots up to 12 crossings---the exponent of $q$ in $f_0$ is determined by the genus $\genusK$ and the \emph{Hopf invariant} $\hopfK$ of the fibered surface. 

The Hopf invariant, introduced by Rudolph as `enhanced Milnor number' \cite{Rudolph1987}, is computed from the contact structure supported by the fibered surface in $S^3$. It detects information beyond classical invariants and, in particular, it can distinguish fibered knots with the same genus and Alexander polynomial. Our result therefore identifies the Hopf invariant as an essential geometric ingredient controlling the leading behavior of $F_K$.
\begin{theorem}\label{thm:nice-L-monomial_intro}
    For a nice knot $K$, the coefficient of the leading term of $F_K$ is a monomial in $q$. In other words, we can write
    \begin{equation}\label{eq:FK-nice-knot_intro}
        F_K = x^\frac12 (\pm q^{\el(K)} x^{\dK-1} + \text{higher order terms in } x),
    \end{equation}
    where $\el(K)$ is an integer-valued invariant of the knot $K$. For a homogeneous braid knot, or a fibered knot of at most 12 crossings,
    \begin{equation}\label{eq:l_hopf_intro}
        \el(K)=\genusK-\hopfK.
    \end{equation}
\end{theorem}

This result was motivated by a recent theorem of L\'opez Neumann and van der Veen in \cite{NeumannVanDerVeen24}, which showed that Hopf invariant appears in the top coefficient of the Akutsu–Deguchi–Ohtsuki (ADO) polynomials \cite{Akutsu1992,Murakami2001} (see \cref{thm:NvdV}). At the same time, ADO polynomials are conjectured to be the radial limits of the $F_K$ series at roots of unity \cite{GukovNakajima2021} (see \cref{conj:FK-ADO}).

Habiro in \cite{Habiro2007unified}  proved the existence of the `cyclotomic expansion' of the colored Jones function:
\begin{equation}\label{eq:Habiro-series-intro}
    J_K(y,q)=\sum_{n=0}^\infty \an{n}{K}\prod_{i=1}^n(y+2-q^i-q^{-i})
\end{equation}
with $\an{n}{K} \in \Z[q,q^{-1}]$. Colored Jones polynomials are obtained as specializations of the form $J_K(q^n+q^{-n}-2,q)$. Similarly, the sequence of \emph{Habiro coefficients} $\an{n}{K}$ is uniquely determined by the sequence of colored Jones polynomials. From \cref{thm:nice-L-monomial_intro}, together with the relation between the expansions of $J_K$ and $F_K$ at $q=1$, we find an explicit formula for the Hopf invariant in terms of the Alexander polynomial and Habiro coefficients, or, equivalently, colored Jones polynomials:

\begin{theorem}\label{thm:Hopf-Habiro_intro} Let $K$ be a homogeneous braid knot, or a fibered knot of at most 12 crossings. Then,
    \begin{equation}\label{eq:hopf_intro}
        \hopfK = \genusK-\sum_{i=0}^{2\genusK-1} \frac{1}{i!}\frac{\partial^i\Delta_K^3(y)}{\partial y^i}\biggr|_{y=0} \frac{\partial \an{2\genusK-i}{K}}{\partial q}\biggr|_{q=1}.
    \end{equation}
\end{theorem}

For example, the Hopf invariant of the knot $8_{10}$ is $4$ and its genus is $3$, while the sum \eqref{eq:hopf_intro} gives
\[
4 = 3-(2030-8505+14076-11508+4662-756).
\]

This result is perhaps surprising. Indeed, L\'opez Neumann and van der Veen in \cite{NeumannVanDerVeen24} write:

\medskip 
\emph{[Their Theorem 1] is the first theorem relating quantum link invariants to fiberedness. No relation seems to be known for the `semisimple' link polynomials, such as colored Jones [\dots].}

\medskip 
Our work shows that the distinction between the semisimple and the nonsemisimple quantum invariants becomes less important when they are all connected via the new quantum invariant $F_K$. 
Naturally, we conjecture the following:

\begin{conjecture}\label{conj:hopf-Jones_intro}
    The formulas \eqref{eq:l_hopf_intro} and \eqref{eq:hopf_intro} hold for all fibered knots.
\end{conjecture}

\begin{figure}
    \centering
    \includegraphics[width=0.175\linewidth, angle=270]{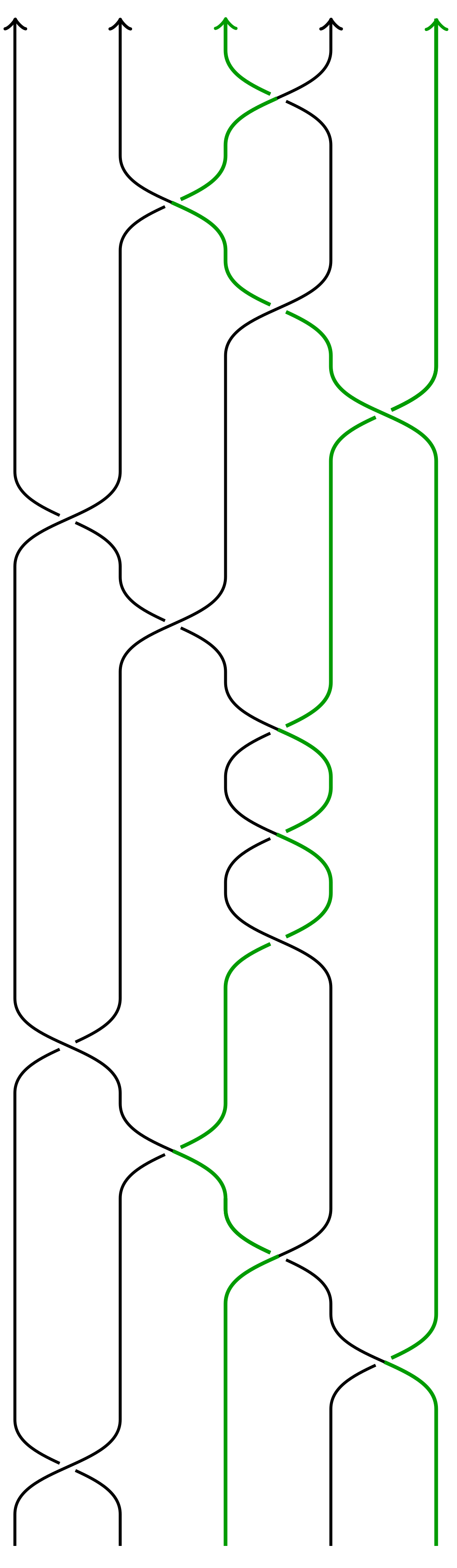}
    \caption{Inversion datum on a braid representative of $12n_{423}$.}
    \label{fig:12n_423_multicycle}
\end{figure}

A knot $K$ is nice if it admits a braid representative equipped with some additional decoration, called \emph{inversion datum}. Roughly, an inversion datum is a choice of a multicycle in the braid allowing certain jumps at each crossing, as in \cref{fig:12n_423_multicycle}. (see \cref{sec:inverted-ss} for definitions). We performed an extensive computer search and found braids with inversion data for all $\num{1246}$ fibered knots up to 12 crossings, extending Park's previous result for knots up to 10 crossings. 

\begin{theorem}\label{thm:nice_12}
    All fibered knots up to 12 crossings are nice.
\end{theorem}

We also confirmed that some fibered knots of crossing number greater than $12$ are nice. Daniel L\'opez Neumann asked us whether the $(2,1)$-cable on the right-handed trefoil is nice, and we find that the answer is positive, which shows that fibered knots that are not obtained solely by plumbing Hopf bands may admit inversion data as well. These facts provide evidence to the `if' part of the following conjecture.

\begin{conjecture}\label{conj:nice-fibered}
    A knot $K$ is nice if and only if it is fibered.
\end{conjecture}

For the `only if' part of the conjecture, we note that no inversion data have been found for non-fibered knots. We expect that the $F_K$ of non-fibered knots is an infinite series in both $x$ and $q$ with integer coefficients \cite{GPV,GPPV,Park20}. However, for a knot to be nice, each $x$-coefficient of $F_K$ must be a Laurent polynomial in $q$. For example, the knot $K=12n_{801}$ is `close' to being fibered, in the sense that its Alexander polynomial is monic with degree equal to the genus of $K$; but the first coefficient $f_0(q)$ of (a natural candidate for) its $F_K$ invariant is the infinite series
\begin{equation}\label{eq:f0_infinite_intro}
f_0(q) =  -q^{-2}(1 + 2 q + q^2 + q^3 + 2 q^4 + 2 q^5 + 2 q^6 + 2 q^7 + \cdots + \num{4112714831} q^{881}+ \cdots).
\end{equation}
Other interesting examples are the fibered knot $12n_{148}$ and the non-fibered knot $13n_{1533}$. While both knots have the same Alexander polynomial, only $12n_{148}$ is nice.

\subsection{Non-fibered knots}

Park constructed the $F_K$ series for a family of twist knots, as well as a handful of other examples, using a `stratified' state sum. We extend the variety of examples and demonstrate two interesting phenomena. 

The first is a relation of the $F_K$ series with the \emph{stability series} of the colored Jones function \cite{GaroufalidisLe2015}. While this relation was previously noted by Gukov and Manolescu for negative torus knots, it appears to be valid in greater generality. We expect that, with a suitable definition of the $F_K$ invariant, $F_K$ will coincide with the stability series (up to a $q$-degree shift by $-\genusK$) for a certain subclass of strongly quasinegative knots. Although it is not clear what is the exact criterion which determines the existence of the stability series, we provide many examples where it exists and coincides with the $F_K$ series.

Secondly, we demonstrate the first examples of knots for which $F_K$ no longer recovers the MMR expansion near $q=1$, one being $12n_{801}$ above. In physics terminology, we interpret this behaviour as the $F_K$ invariant capturing truly non-perturbative effects of the (complex) Chern--Simons theory. 

\subsection{Relation to knot Floer homology}

The knot Floer homology $\widehat{\mathit{HFK}}(K)$ of a knot $K$ is a  doubly graded homology theory categorifying the Alexander polynomial  \cite{Ozsvth2004, RasmussenPhD}. It decomposes as
 \[
 \widehat{\mathit{HFK}}(K) = \bigoplus_{m,n \in \Z} \widehat{\mathit{HFK}}_m(K,n),
 \]
where $n$ denotes the \emph{Alexander grading} and $m$ is the \emph{Maslov (homological) grading}.

Knot Floer homology provides the concordance invariant $\tau(K)$, which bounds the smooth 4-ball genus of the knot $K$ \cite{OzsSza2003}. Moreover, it detects the genus of the knot $\genusK=g$ \cite{OzsSza2008}:
\begin{equation}
    g=\max{\{a\mid\widehat{\mathit{HFK}}(K,a)\neq 0\}}.
\end{equation}
Knot Floer homology also detects fiberedness \cite{Ghiggini2008,Ni2007}:
\begin{equation}
    \dim {\widehat{\mathit{HFK}}(K,g)}=1\iff K \text{ is fibered}.
\end{equation}
The Maslov grading is then given by the Hopf invariant $\lambda(K)=\lambda$ \cite{HFKfibered02}:
\begin{equation}\label{eq:HFK-grading-fibered}
    \widehat{\mathit{HFK}}_m(K,g)\cong\begin{cases}
        \Z \quad &\text{if }m=\lambda, \\
        0 \quad &\text{otherwise.}
    \end{cases}
\end{equation}

For homogeneous braid knots and the first thousand fibered knots, \cref{thm:nice-L-monomial_intro} tells us that $\el$ appears as a (shifted) Maslov grading in knot Floer homology. Motivated to find a connection that may hold in more generality, we study the relation between the GM series and knot Floer homology for non-fibered strongly quasipositive knots.

For the mirrors of these knots, the $F_K$ invariant has as the minimal $q$-exponent of its leading coefficient $-g$. Based on the relation between $\el$ and knot Floer homology for nice knots, and on observation on the knot Floer homology of strongly quasipositive knots, we make the following conjecture:

\begin{conjecture}\label{conj:HFK-SQP_intro} 
Let $K$ be a strongly quasipositive knot with genus $g$. Then the minimal Maslov grading of the top Alexander piece of $\widehat{\mathit{HFK}}$ is zero:
\begin{equation}
    \min{\{m\mid\widehat{\mathit{HFK}}_m(K,g) \neq 0\}} = 0.
\end{equation}
\end{conjecture}

We verified the conjecture for all \num{22009} non-alternating strongly quasipositive knots with at most 16 crossing, tabulated by Stoimenow \cite{StoimenowKnotData}. Standard results on knot Floer homology show that the conjecture holds for fibered or alternating strongly quasipositive knots \cite{HFKalternating03,HFKfibered02}.

As a consequence, we expect that the following holds for fibered or strongly quasinegative knots:
\begin{equation}\label{eq:l_hfk_intro}
\begin{split}
    \el= g-\min{\{m\mid \widehat{\mathit{HFK}}_m(K,g)\neq 0\}}.
\end{split}
\end{equation}
In particular, using knot Floer homology, we are able to provide a potential explanation for the general behavior of $\el$.

Interestingly, the $F_K$ invariant also appears to be sensitive to more than just the minimal Maslov grading in knot Floer homology. In particular, for strongly quasinegative knots, we observed the following pattern: Whenever the Maslov support of $\widehat{\mathit{HFK}}(K,g)$ is 1-dimensional, $F_K$ (or the stability series) recovers the MMR expansion. However, if $\widehat{\mathit{HFK}}(K,g)$ is supported in more than one Maslov grading, $F_K$ may not recover the MMR expansion. We believe that a possible reason behind the similarity between knot Floer homology and $F_K$ could be their behavior under Murasugi sum. While we refer the reader to \cref{rmk:Murasugi} for more details, the study of this phenomenon is beyond the scope of this paper.

\subsection{Slopes}

In \cite{garoufalidis2010slopes}, Garoufalidis proposed a conjecture relating the growth of minimal $q$-powers in the colored Jones function to boundary slopes of essential surfaces in the knot complement. It is natural to expect that a similar statement holds for the $F_K$ invariant.

As the coefficient sequence $f_n(q)$ of $F_K(x,q)$ is $q$-holonomic, the sequence of minimal $q$-degrees of $f_n(q)$ is a quadratic quasi-polynomial for large enough $n$ \cite{Garoufalidis2011}:
\begin{equation}
    \deg_q(f_n) = a(n)n^2+b(n)+c(n) \text{ for } n \gg 0,
\end{equation}
where $a(n)$, $b(n)$ and $c(n)$ are periodic functions. Consider the finite collection of the \emph{$F_K$-slopes} of the knot $K$:
\begin{equation}
    \fs(K) \coloneqq \left\{\frac{1}{a(n)}\right\}_{n \in \N}.
\end{equation}
Let $\bs(K)$ be the collection of boundary slopes of essential surfaces in $S^3 \setminus K$.

\begin{conjecture}[Slope conjecture for $F_K$]\label{conj:slope}
    $F_K$-slopes are boundary slopes; i.e.
    \[
    \fs(K) \subset \bs(K).
    \]
\end{conjecture}

As in the case of the colored Jones function, we observe that the leading coefficient $a(n)$ is a constant, and hence each knot determines a single slope. One may similarly study the values and the common period of $b(n)$ and $c(n)$. However, in contrast with the Jones slopes, $a(n)$ is mostly fractional for small knots, and therefore the period is usually larger than 1. Moreover, the $F_K$ slopes differ substantially from the Jones slopes. As for the Jones slopes, the following problem is quite mysterious, in this case even for alternating knots:
\begin{problem}
Understand which slopes are selected by $F_K$ from the set of all boundary slopes.
\end{problem}

We collect the slopes for fibered Montesinos knots up to 10 crossings in \cref{sec:slopes}. 

\subsection{Conventions}\label{subsec:conventions}

Let $q$ be an indeterminate. For $n \in \Z$ and $k \in \Zpos$, we define
\begin{equation*}
    \sm{n}=q^{\frac n 2}-q^{-\frac n 2}, \enspace [n]=\frac{\sm{n}}{\sm{1}},
    \enspace \qbinom{n}{k}=\frac{[n]\cdots [n-k+1]}{[k]\cdots[1]}, \enspace \qbinom{n}{k}_q =q^{\frac{(n-k)k}{2}}\qbinom{n}{k} = \frac{(q^{n};q^{-1})_k}{(q;q)_k}. 
\end{equation*}
where $(x;q)_k$ is the $q$-Pochhammer symbol $(x;q)_k=(1-x)(1-q x)\cdots (1-q^{k-1}x)$.

Throughout the paper, let $x$ and $y$ be two indeterminates related by $y= x+x^{-1}-2.$ Let $\Delta_K(x) \in \Z[x+x^{-1}]$ be the (symmetrized) Alexander polynomial of a knot $K$, satisfying $\Delta_K(1) = 1$. We abuse notation and write $\Delta_K(y)$ for the same polynomial $\Delta_K(x)$ expressed in terms of the variable $y$. In other words, $\Delta_K(y)$ is the Alexander--Conway polynomial evaluated at $y^{1/2}$. For example, the Alexander polynomial of trefoil knot is
\[
\Delta_{3_1}(x) = x-1+x^{-1}
\]
and we have $\Delta_{3_1}(y) = -1+y$. We use the same convention for other polynomials in $\Z[x+x^{-1}]$.

Given a series $f \in A\lr{x^{1/2}}$, where $A$ is any coefficient ring, we denote by $\lt(f)$ the leading term of $f$, which is the monomial with minimal $x$-power. The \emph{degree} of $f$, denoted by $\dx(f)$, is the minimal $x$-power. For example, $\lt\!\big((1+q^2)x^{-1}+1+q x^{1/2}\big) = (1+q^2)x^{-1}$ and $\dx(\Delta_{3_1}(x)) = -1$. We apply this in particular to polynomials $P(x)\in\Z[x+x^{-1}]$, so that $P(x)\in x^{\dx(P(x))}\Z[x]$. 
The degree $\dx{}$ defines a valuation on the ring $\Q(q)\lr{x}$ of Laurent series whose coefficients are rational functions of $q$, as well as on $\Q(x,q)$ via the embedding 
\begin{equation*}
    \Q(x,q) \hookrightarrow \Q(q)\lr{x}
\end{equation*}
given by the Laurent expansion at $x=0$.

For $n\geq 0$, we denote by $\J{n}{K}$ the $n$-th colored Jones polynomial of $K$, normalized so that $\J{n}{U}=[n+1]$. We follow the conventions of \cite{GM, Habiro2007unified}---the variable $q$ is $A^4$ where $A$ is the Kauffman bracket variable. It is opposite to the convention of most knot theory textbooks and Khovanov homology literature, where $q=A^{-4}$. We also denote by $\Ju{n}{K}$ the unnormalized $n$-th colored Jones polynomial $\J{n}{K}/[n+1]$, such that $\Ju{n}{U}=1$. For right-handed (positive) trefoil $3_1^r$, we have
\begin{align*}
\J{1}{3_1^r} &= q^{-\frac12}+q^{-\frac32} +q^{-\frac52}-q^{-\frac92}\\
&= (q^\frac12+q^{-\frac12})(q^{-1}+ q^{-3}-q^{-4}) = [2]\,\Ju{1}{3_1^r}.
\end{align*}

\subsection*{Acknowledgments} We would like to thank Sergei Gukov and Sunghyuk Park for useful discussions and Davide Passaro for his help with computer calculations. We are grateful to Dror Bar-Natan, Ciprian Manolescu, Daniel L\'opez Neumann and Matthias Storzer for valuable feedback on the first version of this paper. Josef Svoboda was supported by the Simons Collaboration grant New Structures in Low-Dimensional Topology. Lara San Martín Suárez received the support of a fellowship from ``la Caixa'' Foundation (ID 100010434), with fellowship code LCF/BQ/EU23/12010094.
\section{Preliminaries}

\subsection{Inverted state sum}\label{sec:inverted-ss}

\begin{figure}
    \centering
    \includegraphics{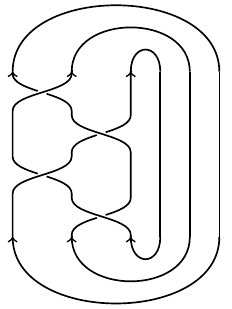}
    \vspace{-5pt} 
\caption{The figure-eight knot $4_1$ represented as a closure of the braid $\beta = \sigma_2^{-1} \sigma_1 \sigma_2^{-1} \sigma_1$.}
\label{fig:braid_closure}
\end{figure}

We recall the inverted state sum and the definition of the $F_K$ invariant from \cite{Park21}.
Let $B_n$ be the braid group on $n$ strands and denote by $\sigma_1,\dots,\sigma_{n-1}$ its usual (Artin) generators. 
Every knot $K$ can be represented as a closure of some braid $\beta$, as depicted in \cref{fig:braid_closure}.
For the rest of this section, we fix a braid $\beta$ which represents a knot $K$. 

Denote by $V_\bd$ the set of crossings in $\bd$ and by $E_\bd$ the set of segments in the braid connecting the crossings. By taking the closure of the braid, we identify the very top segments with the very bottom ones (\cref{fig:braid_state_sum}), such that every segment in $E_\bd$ is associated to exactly two crossings in $V_\bd$. In \cref{fig:crossing-graph} we represent an element of $V_\bd$ as a vertex and the associated segments in $E_\bd$ as directed edges.

\subsubsection{Inversion data}

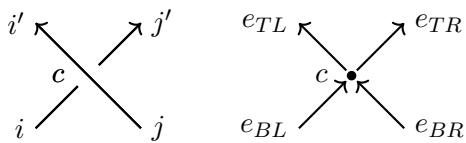
\begin{figure}
    \centering
      \begin{tikzpicture}[scale=0.7,line width=.8pt]
        
        \draw (-1,-1) -- (-0.2,-0.2);
        \draw[->] (0.2,0.2) -- (1,1);
        
        \draw[<-] (-1,1) -- (1,-1);
    
        \node[left] at (-0.25,0) {$c$};

    \begin{scope}[shift={(5,0)}]
        
        \draw[->] (-1,-1) -- (-0.1,-0.1);
        \draw[->] (0.1,0.1) -- (1,1);
        \filldraw (0,0) circle (2pt);
        
        \draw[<-] (-1,1) -- (-0.1,0.1);
        \draw[<-] (.1,-.1) -- (1,-1);
        \node[right] at (1,-1) {$e_{BR}$};
        \node[right] at (1,1) {$e_{TR}$};
        \node[left] at (-1,-1) {$e_{BL}$};
        \node[left] at (-1,1) {$e_{TL}$};
        \node[left] at (-0.25,0) {$c$};
    \end{scope}

        \draw (-1,-1) -- (-0.2,-0.2);
        \draw (0.2,0.2) -- (1,1);
        
        \draw (-1,1) -- (1,-1);
    
        \node[left] at (-0.25,0) {$c$};
        \node[right] at (1,-1) {$j$};
        \node[right] at (1,1) {$j'$};
        \node[left] at (-1,-1) {$i$};
        \node[left] at (-1,1) {$i'$};
    \end{tikzpicture}
    \caption{Every crossing $c\in V_\bd$ is connected to four adjacent segments in $E_\bd$, denoted by $e_{BR}$, $e_{BL}$, $e_{TR}$ and $e_{TL}$ (bottom/top, right/left). The components of a state at a crossing $c$ are denoted by $i,j,i',j'$.}
    \label{fig:crossing-graph}
\end{figure}

\begin{definition}
    An \emph{inversion datum} $\id$ on a braid diagram $\bd$ is an assignment of a $+$ or $-$ sign to each segment of $\bd$:
    \[
        \id: E_\bd \to \{+,-\}.
    \]
    It must satisfy the following conditions: At each crossing, the allowed signs in the segments adjacent to the crossing (the images $\id(e_{BR})$, $\id(e_{BL})$, $\id(e_{TR})$ and $\id(e_{TL})$ as in \cref{fig:crossing-graph}) are
    \begin{equation}\label{eq:allowed-signs-crossings}
    \begin{split}
        &\stackanchor{-+}{+-},\stackanchor{+-}{-+},\stackanchor{--}{--},\stackanchor{-+}{-+},\stackanchor{++}{++} \quad \text{for a positive crossing,} \\
        &\stackanchor{-+}{+-},\stackanchor{+-}{-+},\stackanchor{--}{--},\stackanchor{+-}{+-},\stackanchor{++}{++} \quad \text{for a negative crossing.}
    \end{split}
    \end{equation}
    In particular, at every crossing the number of incoming and outgoing adjacent segments with a $-$ sign coincides.
\end{definition}
In formulas involving $\id$, we identify $\pm$ with $\pm1$. 

\begin{definition}\label{def:state-compatible}
    A \emph{state} $s$ is an assignment of an integer to each segment of the braid
    \[
        s: E_\bd \to \Z,
    \]
    which satisfies the conservation property $s(e_{BL})+s(e_{BR})=s(e_{TL})+s(e_{TR})$ at every crossing.
\end{definition}

For a state $s$ and a crossing $c$, it is convenient to denote the values of $s$ on the adjacent edges by $i,i',j,j'$ as in \cref{fig:crossing-graph} (bottom left, top left, bottom right and top right segments). The dependence on the state $s$ is implicit.

\subsubsection{$R$-matrices}

Denote by $\sgn(c) \in \{1,-1\}$ the sign of the crossing $c$ (positive or negative). For a given crossing $c\in V_\bd$, an inversion datum $\id$, and a state $s$, Park associated a rational function, \emph{the extended $R$-matrix} $R^{\sgn(c)}$, dependent on the values of $s$ on the adjacent edges to $c$ (see \cref{tab:R-matrix-formula}). We consider $F_K$ as a power series in $x$, as opposed to $x^{-1}$, which leads to the opposite signs in the $x$-exponent compared to \cite[eq. 2]{Park21}.

\begin{table}
\centering 
\renewcommand{\arraystretch}{1.5}
\begin{tabular}{l c c}  
    extended $R$-matrix & condition & inv. datum \\
    \midrule[0.3pt]
    \multirow{2}{*}{$R^{+1}(s) = \displaystyle \delta_{i+j}^{i'+j'}q^{jj'}(qx)^{\frac{j+j'+1}{2}}  \qbinom{i}{i-j'}_q (q^{j+1} x;q)_{i-j'}$} 
    & $i\geq j'\geq 0$ & \stackanchor{++}{++}, \stackanchor{-+}{+-} \\ 
    & $0 > i \geq j'$ & \stackanchor{--}{--}, \stackanchor{+-}{-+} \\
    \midrule[0.3pt]
    \rule{0pt}{4ex}$\displaystyle R^{+1}(s) =\delta_{i+j}^{i'+j'}q^{jj'}(qx)^{\frac{j+j'+1}{2}}  \qbinom{i}{j'}_q \frac{1}{(q^{j} x;q^{-1})_{j'-i}}$
    & $j' \geq 0 > i$ & \stackanchor{-+}{-+} \\ 
    \midrule[0.3pt]
    \multirow{2}{*}{$\displaystyle R^{-1}(s) = \delta_{i+j}^{i'+j'}q^{-ii'}(qx)^{-\frac{i+i'+1}{2}}  \qbinom{j}{j-i'}_{q^{-1}} (q^{-i-1}x^{-1} ;q^{-1})_{j-i'}$} 
    & $j\geq i'\geq 0$ & \stackanchor{++}{++}, \stackanchor{+-}{-+} \\
    & $0 > j \geq i'$ & \stackanchor{--}{--}, \stackanchor{-+}{+-} \\
    \midrule[0.3pt]
    $\displaystyle R^{-1}(s) =\delta_{i+j}^{i'+j'}q^{-ii'}(qx)^{-\frac{i+i'+1}{2}}  \qbinom{j}{i'}_{q^{-1}} \frac{1}{(q^{-i} x^{-1};q)_{i'-j}}$ 
    & $i' \geq 0 > j$ & \stackanchor{+-}{+-} \\
\end{tabular}
\caption{The formula for extended $R$-matrices depending on the value of the indices $i,j,i',j'$, also denoted by $\left(R^{\pm1}\right)_{i,j}^{i',j'}$. The $R$-matrix is zero if none of the conditions is satisfied.}
\label{tab:R-matrix-formula}
\end{table}

\subsubsection{Valid states}

\begin{definition}\label{def:state-valid} Given an inversion datum $\id$, we say that a state $s$ is \emph{valid} if 
\begin{enumerate}
    \item the sign of each $s(e)$ coincides with $\id(e)$ for every $e\in E_\bd$, with the convention that the sign of $0$ is $+$,
    \item there is no crossing $c\in V_\bd$ at which ${R^{\sgn(c)}}(s)=0$.
\end{enumerate}
\end{definition}

Of particular importance is the unique valid state only consisting of $0$ and $-1$, which we call the \emph{ground state} and denote by $s_0$. Explicitly, for any $e\in E_\bd$, we have
\begin{equation}\label{eq:ground-state-id}
    s_0(e)\coloneqq \frac{-1+\id(e)}{2}\in\{0,-1\}.
\end{equation}
Clearly, for any valid state $s$, $\abs{s(e)} \geq \abs{s_0(e)}$ for every $e\in E_\bd$.

Let $b_1,b_2,\dots,b_n$ denote the segments associated to the very bottom (and very top) of the braid, as in \cref{fig:braid_state_sum}. We denote by $\Omega(\id)$ the set of all valid states that satisfy $s(b_1)=s_0(b_1)$. 

\subsubsection{State sum}

\begin{figure}
\centering
\begin{tabular}{@{}c c c@{}}
\begin{minipage}[3cm]{0.30\textwidth}
    \centering
    \begin{tikzpicture}
    \pic[braid/number of strands=3, thick] (fig8) at (0,0) 
    {braid =
        {s_1^{-1} s_2 s_1^{-1} s_2}
    };
    
    \path (fig8-1-s) node[left] {$b_1$} to[out=90, in=90, distance=1.5 cm] (3.5,0); 
    
    \draw[thick] (fig8-2-s) node[left] {$b_2$} to[out=90, in=90, distance=1 cm] (3,0); 

    \draw[thick] (fig8-3-s) node[left] {$b_3$} to[out=90, in=90, distance=0.5cm] (2.5,0); 

    \draw[thick] (2.5,0) --(2.5,-4.5);

    \draw[thick] (3,0) --(3,-4.5);

    \path (fig8-rev-1-e) node[left] {$b_1$} to[out=-90, in=-90, distance=1.5 cm] (3.5,-4.5); 

    \draw[thick] (fig8-rev-2-e) node[left] {$b_2$} to[out=-90, in=-90,  distance=1 cm] (3,-4.5); 
    
    \draw[thick] (fig8-rev-3-e) node[left] {$b_3$} to[out=-90, in=-90, distance=0.5cm] (2.5,-4.5);

    \draw[white] (fig8-1-s) to[out=90, in=90, distance=1.5 cm] (3.5,0); 
    \draw[white] (fig8-rev-1-e) to[out=-90, in=-90, distance=1.5 cm] (3.5,-4.5); 

    \draw[thick] (fig8-rev-1-e) -- ++(0,-0.3);
    \draw[<-] (fig8-rev-1-e) -- ++(0,-0.3);
    \draw[<-] (fig8-rev-2-e) -- ++(0,0);
    \draw[<-] (fig8-rev-3-e) -- ++(0,0);

    \draw[<-] (fig8-rev-1-s) -- ++(0,-0.3);
    \draw[<-] (fig8-rev-2-s) -- ++(0,0);
    \draw[<-] (fig8-rev-3-s) -- ++(0,0);
\end{tikzpicture}
    \vfill
    \vspace{-20pt}
    \caption*{(a)}
\end{minipage}
&
\begin{minipage}[3cm]{0.30\textwidth}
    \centering
    \begin{tikzpicture}
    \pic[braid/number of strands=3, thick] (fig8_2) at (0,0) 
    {braid =
        {s_1^{-1} s_2 s_1^{-1} s_2}
    };

    \node[left] at (fig8_2-2-0) {${+}$};
    \node[left] at (fig8_2-1-0) {${+}$};
    \node[left] at (fig8_2-2-1) {${+}$};
    \node[left] at (fig8_2-2-2) {${+}$};
    \node[left] at (fig8_2-2-3) {${+}$};
    \node[left] at (fig8_2-1-3) {${-}$};
    \node[left] at (fig8_2-3-3) {${+}$};
    \node[left] at (fig8_2-1-1) {${+}$};
    \node[left] at (fig8_2-1-2) {${-}$};
    \node[left] at (fig8_2-3-1) {${-}$};
    \node[left] at (fig8_2-1-4) {${+}$};
    \node[left] at (fig8_2-2-4) {${-}$};
    \node[left] at (fig8_2-3-2) {${+}$};
    
    \path (fig8_2-1-s) to[out=90, in=90, distance=1.5 cm] (3.5,0); 


    \draw[thick] (fig8_2-2-s) to[out=90, in=90, distance=1 cm] (3,0); 

    \draw[thick] (fig8_2-3-s) to[out=90, in=90, distance=0.5cm] (2.5,0); 


    \draw[thick] (2.5,0) --(2.5,-4.5);

    \draw[thick] (3,0) --(3,-4.5);



    \draw[thick] (fig8_2-rev-2-e)  to[out=-90, in=-90,  distance=1 cm] (3,-4.5); 
    
    \draw[thick] (fig8_2-rev-3-e)  to[out=-90, in=-90, distance=0.5cm] (2.5,-4.5);

    \draw[white] (fig8_2-1-s) to[out=90, in=90, distance=1.5 cm] (3.5,0); 
    \draw[white] (fig8_2-rev-1-e) to[out=-90, in=-90, distance=1.5 cm] (3.5,-4.5); 

    \draw[<-] (fig8_2-rev-1-e) -- ++(0,-0.3);
    \draw[<-] (fig8_2-rev-2-e) -- ++(0,0);
    \draw[<-] (fig8_2-rev-3-e) -- ++(0,0);

    \draw[<-] (fig8_2-rev-1-s) -- ++(0,-0.3);
    \draw[<-] (fig8_2-rev-2-s) -- ++(0,0);
    \draw[<-] (fig8_2-rev-3-s) -- ++(0,0);
\end{tikzpicture}
    \vfill
    \vspace{-20pt}
    \caption*{(b)}
\end{minipage}
&
\begin{minipage}[3cm]{0.30\textwidth}
    \centering
    \includegraphics{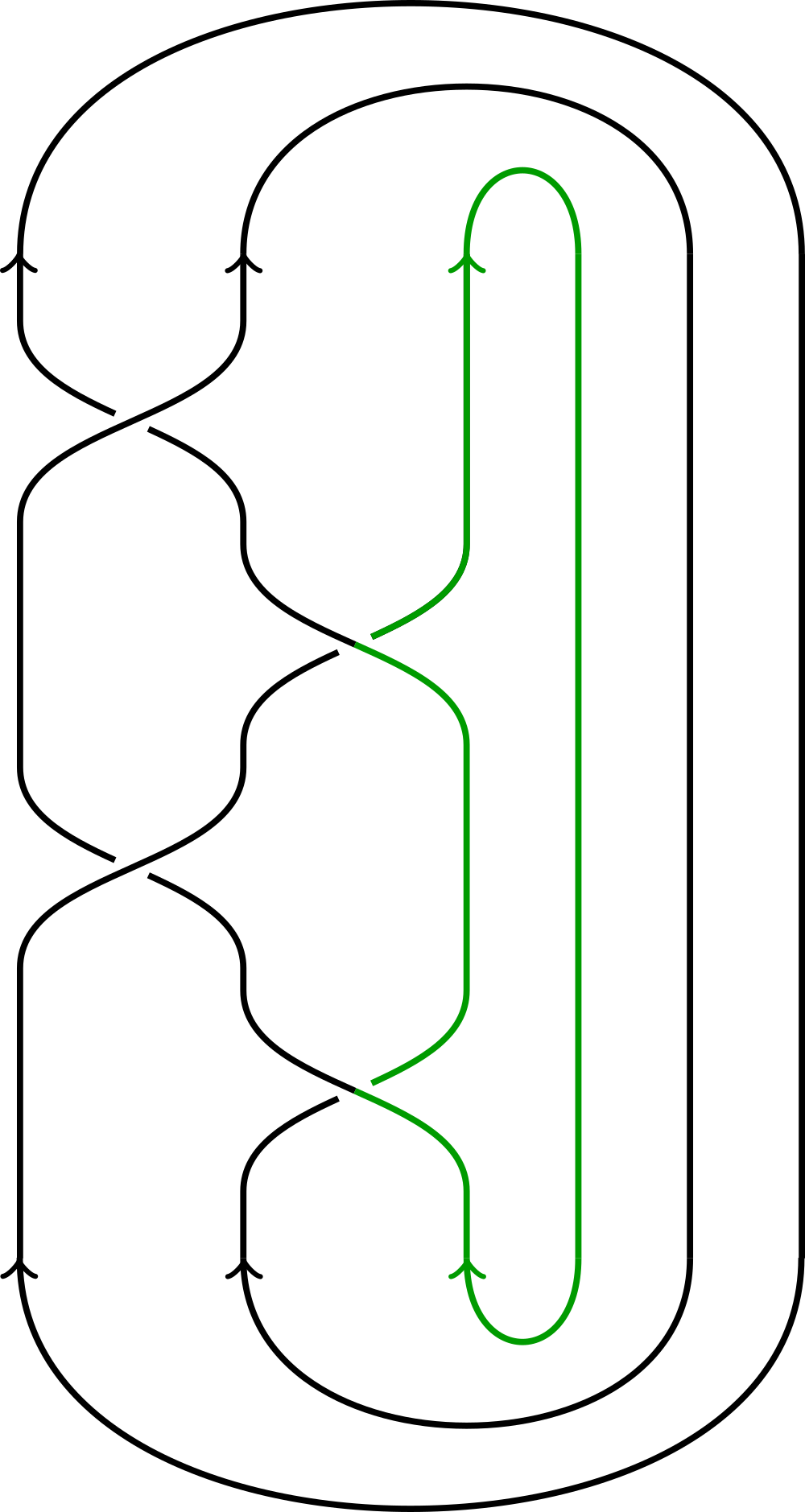}
    \vfill
    \caption*{(c)}
\end{minipage}
\end{tabular}

\caption{(a) The segments $b_1, \dots, b_n$ at the bottom (and top) of the braid diagram $\bd$. 
(b) Inversion datum associated to a braid. 
(c) Each inversion datum can be encoded as an oriented simple multicycle in the braid (in green), allowing jumps from the top to the bottom strands.}
\label{fig:braid_state_sum}
\end{figure}

The product of the extended $R$-matrices over $V_\bd$ associates to any valid state $s\in\Omega(\id)$ a rational function of $x$ and $q$:
\begin{equation}
    P\colon \Omega(\id) \to \Q(x,q).
\end{equation}
defined via
\begin{equation}\label{eq:product}
    P(s)= \left(\prod_{j=2}^{n} x^{-\frac12} q^{-\frac{1}{2}-s(b_j)}\right) \left(\prod_{c\in V_\bd} R^{\sgn(c)}(s)\right).
\end{equation}
The first product arises from taking the quantum trace over the closed strands. The second product consists of the matrix elements of the map we are tracing over. As is customary in such state sum models, we leave the leftmost strand open (see \cref{fig:braid_state_sum}), and for this reason the first product starts at $j=2$.

\begin{definition}[\cite{ParkThesis}]\label{def:inverted-state-sum}
Given an inversion datum $\id$, \emph{the inverted state sum} is the sum
\begin{equation}\label{eq:inverted-state-sum-def}
    Z^{\text{inv}}(\bd) = (-1)^{\mathbf{s}} \sum_{s \in \Omega_1(\id)} P(s).
\end{equation}
Each summand is understood as the Laurent expansion of the rational function $P(s)$ at $x=0$. The overall sign $(-1)^{\mathbf{s}}$ is the sign of the leading coefficient of the Alexander polynomial, see also \cref{rmk:multicycle} or \cite[Thm.~2]{Park21}. We sometimes omit this overall sign as it is not essential for our discussion.
\end{definition}

The minimal $x$-degrees of the extended $R$-matrices are 
\begin{equation}\label{eq:min_x_degree_R}
    \dx((R^{\pm 1})_{i,j}^{i',j'}) = \pm\frac{j+j'+1}{2}
\end{equation}
Therefore, the minimal $x$-degree of $P(s)$ is
\begin{equation}\label{eq:min_x_degree_P}
    \dx(P(s)) = -\frac{n-1}{2}+\sum_{c \in V_\bd} \sgn(c)\frac{j+j'+1}{2} \in\Z.
\end{equation}

In what follows, we often think of $\dx(P(s))$ as an affine function in the variables $s(e)$ for $e \in E_\bd$. Based on degree considerations, Park defined a class of knots for which $Z^{\text{inv}}(\bd)$ does not depend on the order of the summation, and gives a well-defined element of $\Z[q^{\pm}][[x]]$.

\begin{definition}\label{def:nice-knot}
A knot $K$ is \emph{nice} if there is a braid representative $\bd$ and an inversion datum $\id$, such that the following are satisfied:
\begin{enumerate}
    \item The minimal $x$-degrees of all valid states are nonnegative:
    \[ \dx(P(s)) \geq 0, \quad s \in \Omega(\id). \]
    \item Only finitely many states  contribute to a given $x$-degree in the inverted state sum\footnote{Park calls the inverted state sum `absolute convergent' in this case.}. In other words, for every $M \geq 0$, we have that
    \begin{equation*}
    \lvert\{s \in \Omega(\id) \mid \dx(P(s))\leq M \}\rvert< \infty.
    \end{equation*}
\end{enumerate}

\end{definition}
For a nice knot $K$, the $F_K$ invariant is defined as
\begin{equation}\label{eq:FK-state-sum}
    F_K(x,q) \coloneqq (x^{1/2}-x^{-1/2}) Z^{\text{inv}}(\bd).
\end{equation}

We also write $\Fu_K(x,q) \coloneqq Z^{\text{inv}}(\bd)$ (the `unnormalized' $F_K$ invariant). The normalization of $F_K$ is useful in surgery formulas \cite{GM}, while $\Fu_K$ is best suitable for the comparison with the MMR expansion (see \cref{sec:MMR}).

\begin{remark} \cref{def:nice-knot} may be extended to links with more than one component. In this paper we only consider knots for simplicity.
\end{remark}

\begin{remark}\label{rmk:multicycle} Every inversion datum can be interpreted as an oriented, simple multicycle in a braid, where at each crossing we allow `jumps' from the incoming top segment to the outgoing bottom segment (see \cref{fig:braid_state_sum,fig:12n_423_multicycle} for examples and \cite{Park21} for further details). From this point of view, the integer $\mathbf{s}$ in \eqref{eq:inverted-state-sum-def} is the number of closed components of the associated multicycle. This perspective plays a key role in Park's proof of \cref{thm:FK-nice-knot-MMR}, and in particular explains why, at a positive crossing, the signs $\stackanchor{-+}{-+}$ are allowed whereas $\stackanchor{+-}{+-}$ are not.
\end{remark}

\section{The leading term}

\subsection{The leading term and the ground state}

We write the $F_K$ series in the form
\begin{equation}\label{eq:F_K-f0}
    F_K(x,q) = x^{\frac12} \left(f_0(q) x^{d_K-1} + f_1(q) x^{d_K} + \cdots \right), \quad f_0(q) \neq 0. 
\end{equation} 

By the \emph{leading term} of $F_K$, we mean the term $f_0(q) x^{\dK-\frac12}$ associated with the smallest $x$-degree $\dK-\frac12$, where $\dK\in \N$ is a nonzero integer. The following proposition describes the states that contribute to $f_0(q)$:

\begin{proposition}\label{prop:ground-state} Let $K$ be a nice knot with a choice of braid diagram $\beta$ and inversion datum $\id$. Then, the ground state is the only state that contributes to the leading term of the inverted state sum.
\end{proposition}

The proof of this proposition involves applying techniques of linear optimization theory, namely Farkas' lemma, which states the following.
 
\begin{lemma}[{\cite[Lem. 1, p. 318]{Gale1951}}]\label{lemma:farkas-original}
    Let $C\in M_{n\times m}(\R)$ and $\vb{b}\in\R^m$. Then, exactly one of the following two assertions is true: 
    \begin{enumerate}
    
        \item There exists an $\vb{x}\in\R^n$ such that $C^T\vb{x}=\vb{b}$ and $\vb{x}\geq 0$.
        \item There exists a $\vb{y}\in\R^m$ such that $C\vb{y}\geq 0$ and $\vb{b}^T\vb{y}<0$.
    \end{enumerate}
\end{lemma}

For our purposes, we need the following integral modification of the lemma:

\begin{lemma}\label{lemma:farkas}
Let $f:\Z^n\rightarrow \Z$ be an affine function. Let $s \in \{1,2,\dots,n\}$ and let $L$ be a collection of ordered tuples $L\subset [1,n]\times[1,n]$. Let $\mathcal{S}$ be the set of $(a_1,\dots,a_n) \in \Z^n$ satisfying
\[
a_1,\dots,a_s\geq 0, \quad a_{s+1},\dots,a_n\leq0, \quad a_j-a_k\geq 0 \text{ for every } (j,k)\in L.
\]
Assume that $f(a_1,\dots,a_n)\geq 0$ for every $(a_1,\dots,a_n)\in\mathcal{S}$. Then, there exist $\alpha_i,\gamma_{j,k}\in\R_{\geq 0}$ and $\beta\in\Zpos$ such that 
\begin{equation}
    f(a_1,\dots,a_n)=\sum_{i=1}^n \alpha_i \abs{a_i}+\sum_{(j,k) \in L} \gamma_{j,k} (a_j-a_k) + \beta
\end{equation}
for all $(a_1,\dots,a_n)\in\mathcal{S}$.
\end{lemma}

\begin{proof} Assume first that $f$ is linear ($f(0,\dots,0)=\beta=0$) and write $f(a_1,\dots,a_n)=\vb{b}^T\vb{a}=\sum_{i=1}^n b_i a_i$, for some vector $\vb{b} \in \R^n$. 
The set $\mathcal{S}$ consists of vectors $\vb{a}=(a_1,\dots,a_n) \in \Z^n$ that satisfy $C  \vb{a} \geq 0$, where $C$ is the matrix of size $n \times (n + \abs{L})$ of the form:
\begin{equation*}
    C=\begin{pmatrix} 
    1 & &  &  &  &  \\
     & \ddots  &  &  &  &  \\
     & & 1 &  &  &  \\
    &  &  & -1 & &  \\
     & & & & \ddots &  \\
     &  &  &  &  &  -1 \\
    1 & \dots & -1 & & \dots &  \\
    & & & \dots & & \\
\end{pmatrix}
\end{equation*}
The first $n$ rows of $C$ form a diagonal matrix, with diagonal value $+1$ for the first $s$ rows and $-1$ for the rest. The other rows are given by tuples in $L$ in an obvious way. As $C$ contains a diagonal submatrix of size $n$, $C \vb{a}=0$ for $\vb{a}\in\R^n$ if and only if $\vb{a}=\vb{0}$.

The condition $f(a_1,\dots,a_n)\geq 0$ for every $(a_1,\dots,a_n)\in\mathcal{S}$ is satisfied if and only if there is no $\vb{a}\in \Z$ satisfying both $C \vb{a}\geq 0$ and $\vb{b}^T\vb{a}< 0$. In order to apply \cref{lemma:farkas-original}, we need to show that there is no such $\vb{a}$ in $\R^n$. 

First, any such $\vb{a}\in\Q^n$ would give an integral solution of the inequalities by clearing the denominators.
Now assume that there is an $\vb{a}\in\R^n$ satisfying both $C \vb{a}\geq 0$ and $\vb{b}^T\vb{a}< 0$. Actually, we must have that $C \vb{a}> 0$, since $C \vb{a}= 0$ would imply that $\vb{a}=\vb{0}\in\Z^n$, which is a contradiction. Therefore, there is an $\vb{a}\in\R^n$ satisfying both $C \vb{a}> 0$ and $\vb{b}^T\vb{a}< 0$. As this is an open condition, we can find some $\vb{y} \in \Q^n$ satisfying it, obtaining again a contradiction.

We can now apply Farkas' lemma. By \cref{lemma:farkas-original} there exist $\lambda_i\geq 0$, ${i=1,\dots, n+\abs{L}}$ such that $\vb{b}=\sum_{i=1}^{n+\abs{L}} \lambda_i \vb{c}_i$, where $\vb{c}_i$ are the rows of $C$. The result follows. If $f(\vb{0}) \neq 0$, we can follow a similar argument with Farkas' lemma in its affine form, see e.g. \cite{notes_farkas}.
\end{proof}

\begin{proof}[Proof of \cref{prop:ground-state}]
    Enumerate the edges in $E_\bd$ by $e_1,\dots,e_{m}$ with $m = \abs{E_\bd}$. By \eqref{eq:min_x_degree_P}, $\dx{P(s)}$ can be thought as an affine function in variables $s(e_1),\dots,s(e_{m})$. We want to apply \cref{lemma:farkas} on $\dx(P(s))$ to deduce that the ground state $s_0$ minimizes $\dx(P(s))$ over the valid states, and it is a unique such state.
    
    For any valid state $s$, the values $s(e)$ for $e\in E_\bd$ satisfy the following relations: For every crossing $c$ such that $\id(e_{BR})=\id(e_{TL})$ and $\id(e_{BL})=\id(e_{TR})$,
    \begin{equation}
        \sgn(c)\left(s(e_{BR})-s(e_{TL})\right)\geq 0 \quad \text{and} \quad \sgn(c)\left(s(e_{BL})-s(e_{TR})\right)\geq 0.
    \end{equation}
    or in the notation $i,j,i',j'$ in \cref{fig:crossing-graph},
    \begin{equation}
        \sgn(c)\left(i-j'\right)\geq 0 \quad \text{and} \quad \sgn(c)\left(i'-j\right)\geq 0.
    \end{equation}

    Denote by $L$ the set of ordered tuples $(r,t)\in\{1,\dots,m\}^2$ such that $s(e_r)-s(e_t)\geq 0$ for every $s\in\Omega(\id)$ as described above.
    Let $\mathcal{S}=\{(s(e_1),\dots,s(e_{m}))\mid s\in\Omega(\id)\}\subset \Z^n$. Setting $a_k = s(e_k)$ for $k=1,\dots m$, $\mathcal{S}$ consists of all vectors $(a_1,\dots,a_{m}) \in \Z^n$ satisfying
    \begin{itemize}
        \item $a_k\geq0$ if $\id(e_k)=+1$,
        \item $a_k<0$ if $\id(e_k)=-1$,
        \item $a_r-a_t\geq 0$ for $(r,t)\in L$.
    \end{itemize}
    
    By the first condition in \cref{def:nice-knot}, the minimal $x$-degree $\dx(P(s))$ of any valid state $s$ is always non-negative. If we view $\dx(P)$ as a function on $a_k=s(e_k)$ for $k=1,\dots,m$, instead of as a function on $s\in\Omega(\id)$, the non-negative condition is equivalent to $\dx(P)\geq 0$ in $\mathcal{S}\subset\Z^n$.
    By \cref{lemma:farkas}, there exist $\alpha_k\in \R_{\geq 0}$, $\gamma_{r,t}\in \R_{\geq 0}$ and $\beta\in\Zpos$ such that $\dx(P(s))$ can be expressed as 
    \begin{equation}
        \dx(P(s))=\sum_{k=1}^{m} \alpha_{k} |s(e_k)| + \sum_{(r,t)\in L}^{}\gamma_{r,t} |s(e_r)-s(e_t)| + \beta.
    \end{equation}
    
    Therefore, for every edge $e_k\in E_\bd$ with $\alpha_k>0$, the value of $s(e_k)$ minimizing $\dx(P(s))$ agrees with the value assigned by the ground state \eqref{eq:ground-state-id}.
    Likewise, for every pair of edges $\{e_r,e_t\}$ such that $\gamma_{r,t}>0$, we must have that $s(e_r)=s(e_t)$ (recall that $(r,t)\in L$ only if $\id(e_r)=\id(e_t)$).
    In particular, the ground state $s_0$ minimizes $\dx(P(s))$ over the valid states as desired.
    
    Assume now for contradiction, that another state $s \neq s_0$ minimizes $\dx(P(s))$.
    Since the values of the state at the edges $e_k$ with $\alpha_k>0$ are uniquely determined, the state $s$ must have $|s(e_k)| > \frac{1-\id(e_k)}{2}$ for some edge $e_k$ with $\alpha_k=0$. The state $s$ cannot be valid: If it was, we could change the value of $s(e_k)$ to an arbitrary number, obtaining infinitely many states with $\dx(P(s))=\dx(P(s_0))$. This would contradict the second condition in \cref{def:nice-knot}. Therefore, $s\notin \Omega(\id)$, and we are done.
\end{proof}

\begin{theorem}\label{thm:nice-L-monomial}
    For a nice knot $K$, the coefficient of the leading term of $F_K$ is a monomial in $q$. In other words, we can write
    \[
    F_K = x^\frac12 \left((-1)^{\mathbf{s}+1}q^{\el(K)} x^{\dK-1} + \text{higher order terms in } x\right),
    \]
    where $\el(K)$ is an integer-valued invariant of the knot $K$.
\end{theorem}

\begin{corollary}\label{cor:FK-min-xdeg}
    $d_K$ is the half-degree of the Alexander polynomial.
\end{corollary}
\begin{proof}
    Follows from the previous theorem and \cref{thm:FK-nice-knot-MMR}. See \cref{sec:MMR} for the specific details.
\end{proof}

\begin{proof}[Proof of \cref{thm:nice-L-monomial}] By \cref{prop:ground-state}, the ground state $s_0\in \Omega(\id)$ is the unique state contributing to the leading term of the state sum. Moreover, for every crossing $c$, the contribution of $R^{\sgn(c)}(s_0)$ to the leading term is a monomial in $q$ (see \cref{tab:R-matrix-formula,tab:xq-contributions}). Therefore, $P(s_0)$ is a monomial in $q$ and $x$, and so is $F_K$.
\end{proof}

\begin{table}
    \centering
    \begin{minipage}{.46\linewidth}
        \begin{tabular}{c|*{3}{c}}

          \begin{tikzpicture}[scale = .25, rotate = 90, baseline=-2pt]
        \draw (-1,-1) -- (-0.2,-0.2);
        \draw (0.2,0.2) -- (1,1);
        \draw (-1,1) -- (1,-1);
        \end{tikzpicture}

        & \stackanchor{--}{--} & \stackanchor{+-}{-+},\stackanchor{-+}{+-} & \stackanchor{-+}{-+},\stackanchor{++}{++}\\
        \midrule
        $R(s_0)$ & $q^{\frac12}x^{-\frac12}$ & 1 & $q^{\frac12}x^{\frac12}$ \\
        \end{tabular}
        \label{tab:xq-pos}
    \end{minipage}
    \begin{minipage}{.5\linewidth}

        \begin{tabular}{c|*{3}{c}}
        
        \begin{tikzpicture}[scale = .25, baseline=-2pt]
        \draw (-1,-1) -- (-0.2,-0.2);
        \draw (0.2,0.2) -- (1,1);
        \draw (-1,1) -- (1,-1);
        \end{tikzpicture}
        
        & \stackanchor{++}{++} & \stackanchor{+-}{-+},\stackanchor{-+}{+-} & \stackanchor{+-}{+-},\stackanchor{--}{--} \\
        \midrule
        $R^{-1}(s_0)$ & $q^{-\frac12}x^{-\frac12}$ & 1 & $q^{-\frac12}x^{\frac12}$ \\
        \end{tabular}
        \label{tab:xq-neg}
    \end{minipage}   
    \caption{Leading term of $R^{\sgn(c)}$ evaluated at the state $s_0$, in terms of the value of $\id$ at a crossing $c$.}
    \label{tab:xq-contributions}
\end{table}

\subsection{Formula for the leading coefficient}

Define $\Rq\in\left\{\pm\frac{1}{2}\right\}$ as the $q$-exponents summarized in Table \ref{tab:xq-contributions}, which depend on the value of $\id$ on two out of the four segments adjacent to a crossing $c$:
\begin{equation}
    \Rq \coloneqq \sgn(c)\frac{1+\id(e_{BR})\id(e_{TR})}{4}\in\left\{\pm\frac12\right\}.
\end{equation}

Then, in the notation of \cref{sec:inverted-ss}, we can express $\el(K)$ in terms of $\bd$ and $\id$ as
\begin{equation}\label{eq:L-general-formula}
    \el(K)=-\sum_{j=2}^n \frac{\id(b_j)}{2}+\sum_{c\in V_\bd} \Rq.
\end{equation}
Here we used \eqref{eq:ground-state-id} to write
\[
\sum_{j=2}^n \left(-\frac{1}{2}-s_0(b_j)\right)=-\sum_{j=2}^n \frac{\id(b_j)}{2}.
\] 

The formula \eqref{eq:L-general-formula} simplifies in the case of homogeneous braid knots. A braid $\bd\in B_n$ is \emph{homogeneous} if it is given by a word in $\{\sigma_1^{\ep_1},\dots,\sigma_{n-1}^{\ep_{n-1}}\}$ for fixed values of $\ep_i=\pm 1$ for $i=1,\dots,n-1$. A \emph{homogeneous braid knot} is a knot $K$ that admits a homogeneous braid representative. The \emph{writhe} of $\beta$ is $w(\bd) = n_+(\bd)-n_-(\bd)$, the difference between the number of positive and negative crossings in $\beta$.

\begin{proposition}\label{prop:L-homogeneous-braid} Let $K$ be a knot represented by a homogeneous braid $\bd\in B_n$. Then,
\begin{equation*}
    \el(K) = \frac{w(\bd)}{2}-\sum_{i=1}^{n-1}\frac{\ep_i}{2}.
\end{equation*}
\end{proposition}

\begin{proof} A homogeneous braid $\bd$ admits a natural inversion datum $\id$, given as follows \cite{Park21}:
\begin{enumerate}
    \item $\id(e_{BR})=\id(e_{TR})=\sgn(c)$ for both edges adjacent to any crossing $c$ on the right.
    \item $\id$ takes a constant value on the edges adjacent on the left to all crossings corresponding to $\sigma_1$ (`leftmost column'), $e_{TL}$ and $e_{BL}$. Here there are two choices, depending on whether we assign $+$ or $-$ sign to the leftmost column.
\end{enumerate}

For this inversion datum, the signs on the segments around a positive crossing are either $\stackanchor{++}{++}$ or $\stackanchor{-+}{-+}$, and both contribute $+\frac12$ to the $q$-exponent. For a negative crossing, the signs are $\stackanchor{+-}{+-}$ or $\stackanchor{--}{--}$, and contribute $-\frac12$ (see \cref{tab:xq-contributions}). The formula \eqref{eq:L-general-formula} takes the form:
\begin{equation*}
\begin{split}
    \el(K)&=-\sum_{j=2}^n\frac{\id(b_j)}{2}+\sum_{c\in V_\bd} \Rq = -\sum_{i=1}^{n-1}\frac{\varepsilon_i}{2}+\sum_{c\in V_\bd} \frac{\sgn(c)}{2} = -\sum_{i=1}^{n-1}\frac{\varepsilon_i}{2}+ \frac{w(\bd)}{2}. \qedhere
\end{split}
\end{equation*}
\end{proof}

\begin{remark}\label{rmk:L-homogeneous-braid-tau-s} For homogeneous braid knots, $\el(K)=\frac{s(K)}{2}=\tau(K)$, where $s$ denotes the Rasmussen's $s$-invariant and $\tau$ denotes the Ozsváth--Szab\'o's $\tau$-invariant \cite{Lobb10,Kawamura2009}. Note that this is not true for an arbitrary fibered knot. In particular, $\el$ is not a concordance invariant, since  $\abs{\el(10_{153})}=1$ for the slice knot $10_{153}$.
\end{remark}

\begin{example} The $F_K$ series of the figure-eight knot $4_1$, with braid representative $\beta = \sigma_1\sigma_2^{-1}\sigma_1\sigma_2^{-1}$, is 
\[
F_{4_1}(x,q) = x^\frac12 (1 + 2x + (q^{-1} +3 + q) x^2  + x^3(2q^{-2}+2q^{-1} +5 +2 q + 2 q^2) + \cdots )
\]
In this case $d_{4_1}=1$. We have $w(\beta)=0,\varepsilon_1=1$ and $\varepsilon_2=-1$, hence $\el(4_1)=0$.
\end{example}

\section{Nice and fibered knots}\label{sec:nice_fibered}

\subsection{Inversion data of fibered knots}

From the relation between $F_K$ and the Alexander polynomial (see \cref{thm:FK-nice-knot-MMR}), it immediately follows that the Alexander polynomial of a nice knot must be monic, just as in the case for fibered knots. It is therefore natural to compare these two families. 

Park demonstrated that every fibered knot with at most 10 crossings is nice. We extended this result to all fibered knots with at most 12 crossings:

\begin{theorem}
    All fibered knots up to 12 crossings are nice.
\end{theorem}

\begin{proof}

The proof was established by a laborious, computer-assisted search. For many knots, their tabulated braid representations do not admit inversion data, so the main challenge for us was to enumerate a variety of new braid representations to test. We focused on obtaining minimal length braid representatives, given that state sum models grow exponentially with the number of crossings.

Our starting point was the list of the $\num{1246}$ fibered knots up to 12 crossings in the KnotInfo database \cite{knotinfo}, along with one of their minimal length braid representatives. For each braid of length $n$, there are $2n$ segments and thus $2^{2n}$ ways to label its segments with a $+$ or $-$ sign. We tested the subset of these labelings for which the $R$-matrices were non-zero and satisfy the criteria of \cref{tab:R-matrix-formula} at each crossing. These are in one-to-one correspondence with the number of simple multicycles in the braid, when we allow jumps from the overstrand to the understrand (\cref{rmk:multicycle}). Given these constraints, for each braid of length $n$, there are at most $2^n$ ways to label its segments. An exhaustive test of segment labelings for each braid in KnotInfo revealed that some of them admitted an inversion datum, but many did not.

For the braids in KnotInfo that do not admit an inversion datum, we used braid moves (braid group relations and Markov moves) to search for different braid presentations. Our most successful search strategy was to start with a known braid and carry out random sequences of braid moves. Each resulting terminal braid in the sequence was normalized by performing obvious reductions and conjugating to the lexicographically smallest braid word. By iterating this search from new braids, we found an expanding list of minimal length braids for a given knot, and the normalization ensured that each new braid we found was sufficiently distant in the search space from our existing list to avoid redundancy. For fibered knots up to 12 crossings, this approach typically yielded a minimal length braid with an inversion datum after at most $10^6$ iterations, each consisting of $10^3$ random braid moves. 

A handful of knots were resistant to this approach, including $12n_{423}$, which appears to have a minimum length braid with 14 crossings.  To find more braid representatives of $12n_{423}$, we created an exhaustive list of braids of length 14, such that any braid of length 14 would normalize to a braid in our list, a mirror image of a braid in our list, or a shorter braid.  Among our $\sim 10^9$ braids of length 14, we computed with SnapPy \cite{SnapPy} that $\sim 10^3$ of them matched $12n_{423}$ on several invariants, including $\widehat{\mathit{HFK}}$ and the Jones polynomial. SnapPy's manifold identification tool recognized a few hundred of these braids as $12n_{423}$ or its mirror, among which we found multiple braids admitting inversion data (an example shown in \cref{fig:12n_423_multicycle}).

While our search for inversion data involved lots of experimentation and ad-hoc tuning of search parameters, the results were straightforward to test.  Our resulting data set is a list of a braid and an inversion datum for each fibered knot up to 12 crossings, and these are verifiable with a separate computation that takes only a few minutes. The braids with inversion data are attached in a \verb|.tsv| file, where the signs of the inversion data are ordered by following the braid from the bottom end of the open strand to its top end.

\end{proof}

\begin{remark}
We have not found any braid admitting more than one inversion data (up to flipping the value of $\id$ on the open strand). This uniqueness, together with \cref{thm:L-Hopf-homogeneous}, suggests that the inversion data should have a geometric interpretation.
\end{remark}

\begin{remark}

In contrast with the theorem above, every fibered knot up to 12 crossings has a braid which does \emph{not} admit an inversion datum. It appears that such braids can be obtained in a systematic way, as follows. For each knot $K$ and braid representative $\beta$ admitting an inversion datum in our data set, we found a modified braid $\beta'$
of the same knot $K$ with no inversion datum. The modified braid $\beta'$ is obtained by two stabilizations (Reidemeister I twists) followed by adding a canceling pair of the last generator and its inverse.  In other words, $\beta'$ is a connected sum of $\beta$ with an unknot, represented by $\sigma_1 \sigma_2 \sigma_2 \sigma_2^{-1}$. 

\end{remark}

\subsection{Relation to the Hopf invariant}\label{subsec:hopf}

Recently, L\'opez Neumann and van der Veen in \cite{NeumannVanDerVeen24} studied the nonsemisimple quantum invariants of Akutsu--Deguchi--Ohtsuki (ADO). They found that, for fibered links, the Hopf invariant appears in the top coefficient of the ADO polynomials. We formulate their main theorem for knots for simplicity:
\begin{theorem}[{\cite[Thm.\! 1]{NeumannVanDerVeen24}}]\label{thm:NvdV}
Let $K \subset S^3$ be a fibered knot. Then, for any $p \geq 2$, the associated ADO polynomial $N_p(K,x)$ has degree $2g(K)(p-1)$ and its top coefficient is given by 
\begin{equation}
    (-\zeta_p)^{g(K)-\lambda(K)},
\end{equation}
where $\zeta_p=e^{2i \pi/p}$ and $\lambda(K)$ is the Hopf invariant of $K$.
\end{theorem}
On the other hand, the polynomials $N_p(K,x)$ conjecturally arise as radial limits $q \nearrow \zeta_p$ of the $F_K$ series. We formulate the conjecture for nice knots, as some examples suggest that it does not hold in general (see \cref{subsec:radial_limits}):
\begin{conjecture}[{\cite[Conj.\! 3]{GukovNakajima2021}}]\label{conj:FK-ADO} For a nice knot $K$,
\[
    {\frac{F_K(x,q)}{x^{1/2}-x^{-1/2}}}\biggr|_{q=\zeta_p} = \frac{N_p(K,x)}{\Delta_K(x^p)}.
\]
\end{conjecture}

Motivated by the above and by \cref{thm:nice-L-monomial}, we make the following conjecture:

\begin{conjecture}\label{conj:L-hopf}
    Let $K$ be a fibered knot. Then, $K$ is nice and
    \begin{equation}
        \el(K)=\genusK-\hopfK,
    \end{equation}
    where $\hopfK$ is the Hopf invariant of $K$. In other words,
    \begin{equation}
        \lt(F_K) = (-1)^{\hopfK+1} q^{\genusK-\hopfK} x^{\genusK-1/2}.
    \end{equation}
\end{conjecture}

This conjecture would extend \cref{thm:NvdV} from roots of unity to generic $q$ in the unit disc. We have a partial result in this direction:

\begin{theorem}\label{thm:L-Hopf-homogeneous} \cref{conj:L-hopf} is true for all homogeneous braid knots, fibered knots with at most 12 crossings, and the $(2,1)$-cable on the right-handed trefoil.
\end{theorem}

\begin{proof} The overall sign of $F_K$ for a nice knot is opposite to the sign of the leading coefficient of the Alexander polynomial, which is known to be $(-1)^{\hopfK}$ for fibered knots. The homogeneous braid case follows from a direct comparison of \cref{prop:L-homogeneous-braid} and \cite[Corollary 4.15]{NeumannVanDerVeen24}. For fibered knots up to 12 crossings, we computed their $\el(K)$ and compared it with the extremal knot Floer homology, which encodes $g(K)$ and $\lambda(K)$ as the Alexander and Maslov grading, respectively. For the $(2,1)$-cable of the right-handed trefoil, we checked that the braid representative
\begin{equation*}
    \sigma_1 \sigma_2 \sigma_3^{-2} \sigma_2 \sigma_1 \sigma_3 \sigma_2 \sigma_3^{-1} \sigma_2\sigma_1 \sigma_3 \sigma_2\sigma_3^{-1} \sigma_2 \sigma_3 \sigma_2^{-1}
\end{equation*}
with the inversion datum $\id\equiv 1$ is nice. In this case, the inverted state sum computation gives
\begin{equation*}
    F_{K}(x,q)=x^\frac12\left(-q^2 x + q^4 x^2 - q^7 x^3 + q^{11} x^4 + (q^5 - q^{16}) x^5 + (-q^{11} + q^{22}) x^6 + O(x^7)\right).
\end{equation*}
\end{proof}

Note that many fibered knots under 12 crossings are not homogeneous knots, and a fortiori not homogeneous braid knots \cite{Abe10-2}; the first such example is the knot $8_{20}$. Moreover, the $(2,1)$-cable on the right-handed trefoil is nice, which shows that fibered knots that are not obtained solely by plumbing Hopf bands may admit inversion data as well. This, and the presence of the Hopf invariant in the $F_K$ invariant of fibered knots, provides support for \cref{conj:nice-fibered}. It also suggests that the inversion datum itself may have a topological interpretation---it seems to be related to the fiber surface, which cannot in general be obtained from the braid diagram by Seifert's algorithm.

\section{Expansions at roots of unity}\label{sec:MMR}

\subsection{Expansion of \texorpdfstring{$F_K$}{F} at 1}

In this section, we implicitly identify rational functions in $x$ with their Laurent expansions at $x = 0$, and similarly Laurent polynomials in $q$ with their power series in $h := q-1$ at $h=0$. We recall the relation between the $F_K$ invariant and the Melvin--Morton--Rozansky (MMR) expansion of a knot $K$.

\begin{theorem}[{\cite[Thm.\! 2]{Park21}}]\label{thm:FK-nice-knot-MMR} Let $K$ be a nice knot. Then, the expansion of the (unnormalized) $F_K(x,q)$ series at $q=1$ coincides with the Melvin--Morton--Rozansky expansion of the colored Jones polynomial at $x = 0$:
\begin{equation}\label{eq:FK-MMR}
    \Fu_K(x,1+h)= \frac{1}{\Delta_K(x)}+\sum_{i=1}^\infty\frac{P^{(j)}_K(x)}{\Delta^{2j+1}_K(x)}h^j.
\end{equation}
\end{theorem}

Both the Alexander polynomial $\Delta_K(x)$ and the polynomials $P^{(j)}_K(x)$ are elements of $\Z[x+x^{-1}]$ (see \cite[Section 7.2]{Habiro2007unified}). Using \cref{thm:nice-L-monomial}, we can derive conditions on the degree and leading coefficient of $P^{(j)}_K$ in terms of $\el(K)$ for a nice knot $K$. Recall our convention that the degree of $P^{(j)}_K(x)$ refers to its minimal $x$-exponent.

\begin{proposition}\label{prop:deg-Pj-MMR-posL} Let $K$ be a nice knot with $\el(K)\geq 0$. For $0 \leq j \leq \el(K)$, we have $\deg_{x} P^{(j)}_K(x) = -2\dK j$  
with leading coefficient $\binom{\el(K)}{j}$. For $j>\el(K)$, we have $
\deg_{x} P^{(j)}_K(x) > -2\dK j$.
\end{proposition}
\begin{proof} 
According to \cref{thm:nice-L-monomial}, the leading term of the series $F_K(x,q)$ is of the form
$\lt(\Fu_K(x,q))=(-1)^{\mathbf{s}} q^{\el(K)}x^{\dK}$,
which expanded at $h=q-1$ near 0 gives
\begin{equation}\label{eq:integral-FK-positive-L}
    \begin{split}
        \lt(\Fu_K(x,1+h))&=(-1)^{\mathbf{s}} (1+h)^{\el(K)} x^{\dK} = (-1)^{\mathbf{s}} \sum_{k=0}^{\el(K)} \binom{\el(K)}{k} h^kx^{\dK}.
    \end{split} 
\end{equation}
Comparing this with the MMR expansion \eqref{eq:FK-MMR} when $h\equiv 0$, we deduce that $\dK$ must be the half-degree of Alexander polynomial of $K$, which shows \cref{cor:FK-min-xdeg}. Moreover, using $\lt(\Delta_K^{2j+1}(x)) = (-1)^{\mathbf{s}} x^{-(2j+1)\dK}$,
we obtain
\begin{equation*}
    \lt\left(\frac{P^{(j)}_K(x)}{\Delta_K^{2j+1}(x)}\right) = (-1)^{\mathbf{s}} x^{(2j+1)\dK}\lt(P_K^{(j)}(x)) = (-1)^{\mathbf{s}} \binom{\el(K)}{j} x^{\dK},
\end{equation*}
which must agree with the $h^j$ coefficient of \eqref{eq:integral-FK-positive-L}. 
\end{proof}

We can do a similar analysis when $\el(K)<0$. In this case, all the polynomials on the right-hand side of \eqref{eq:FK-MMR} have the minimal possible $x$-degree.
\begin{proposition}\label{prop:deg-Pj-MMR-negL} Let $K$ be a nice knot with $\el(K)<0$. Then $\deg_x P_K^{(j)}(x)=-2\dK j$ for every $j\geq0$, 
with leading coefficient $(-1)^j\binom{-\el(K)+j-1}{j}$.
\end{proposition}

\begin{proof} Up to overall sign, the expansion of $\Fu_K(x,q)$ at $h=q-1=0$ reads
\begin{equation}
    \begin{split}
        \lt(\Fu_K(x,1+h))&=(1+h)^{-|\el(K)|} x^{\dK} = \sum_{k=0}^{\infty} (-1)^k\binom{-\el(K)+k-1}{k} h^kx^{\dK}.
    \end{split}
\end{equation}
The rest of the proof carries similarly to the proof of \cref{prop:deg-Pj-MMR-posL}.
\end{proof}

As a special case of Propositions \ref{prop:deg-Pj-MMR-posL} and \ref{prop:deg-Pj-MMR-negL}, we obtain:

\begin{theorem}\label{thm:L-P1} Let $K$ be a nice knot. Then $\el(K)$ is the coefficient of $x^{2\dK}$ in $P^{(1)}_K(x)$. 
\end{theorem}

Note that the bound $\deg P^{(1)}_K\geq -2\genusK$ holds for all knots \cite{Ohtsuki07}. \cref{thm:L-P1} provides a refinement to this bound for nice knots, by showing exactly when it is sharp in terms of $\el(K)$.

\begin{remark}  

Recently, the authors of \cite{BNVdV18} defined a polynomial time invariant $\rho_1(K,x)$, conjecturally equivalent to $P_K^{(1)}(x)$. In \cite[Theorem 54]{BNVdV21} they proved that $\deg \rho_1(K,x)\geq -2\genusK$, by relating $\rho_1$ to the Seifert surface of the knot. If all fibered knots are nice, \cref{conj:L-hopf} provides a geometric interpretation for the degree bound of $P^{(j)}_K$, and conjecturally $\rho_1$, in terms of the genus and the Hopf invariant of the fibered surface. We verified that this is the case for homogeneous braid knots, as well as all fibered knots up to 12 crossings. It would be interesting to study how the Hopf invariant may arise in the invariant $\rho_1$ using Seifert surface arguments.
\end{remark}

\begin{remark} Following a suggestion of Bar-Natan, we looked at the invariant $\theta$ defined in \cite{BNvdV25}. This $\sl_3$ analogue  of $\rho_1$ is expected to agree with the numerator of the $h^1$ term in the $\sl_3$ MMR expansion (cf. \cite[Thm.\! 1.1]{GruSan2025}). \cite[Conj.\! 19]{BNvdV25} predicts that, for a fibered knot $K$, the coefficient of $T_2^{2g}$ in $\theta(T_1,T_2)$ is a shifted Alexander polynomial multiplied by a certain integer $s(K)$. Motivated by \cref{thm:L-P1}, we verified that \[s(K)=\el(K)=g(K)-\lambda(K)\] for all fibered knots up to 12 crossings. Naturally, we expect this to hold for all fibered knots.
\end{remark}

\subsection{Hopf invariant from colored Jones polynomials}\label{subsec:habiro_hopf}

An interesting corollary of \cref{thm:L-P1} is a formula for $\el(K)$ in terms of the Alexander polynomial and the colored Jones polynomials of $K$. In order to formulate it, we recall the \emph{Habiro series} (\cite[Thm.\! 4.5]{Habiro2007unified}) of the knot $K$. It is a series of the form
\begin{equation}
    J_K(y,q)=\sum_{n=0}^\infty \an{n}{K}\prod_{i=1}^n(y+2-q^i-q^{-i}).
\end{equation}
The \emph{Habiro coefficients} $\an{n}{K}\in\Z[q,q^{-1}]$ of $K$ are uniquely determined by the colored Jones polynomials.
\begin{equation}\label{eq:Habiro-Jones}
    \an{n}{K} = \frac{1}{\{2n+1\}_{2n}}\sum_{i=0}^n (-1)^{n-i} \frac{[2i+2]}{[n+i+2]} \qbinom{2n+1}{n+1+i} \J{n}{K}.
\end{equation}
Likewise, the colored Jones polynomials are specializations of the Habiro series.
\begin{equation}
    \J{n}{K} =[n+1] J_K(q^n+q^{-n}-2,q).
\end{equation}

The Habiro series maps to Rozansky’s integral version of the Melvin--Morton expansion of the Jones polynomial, via a natural homomorphism (see \cite[Section 7.2]{Habiro2007unified}). Thus the MMR expansion appears both in the expansions of the Habiro series and the $F_K$ series, after identifying $y$ with $x+x^{-1}-2$, as illustrated in the following diagram:

\begin{center}
    \centering
\begin{tikzpicture}[
  >=Stealth,
  every node/.style={font=\normalsize, align=center}
]

\matrix (m) [matrix of math nodes, row sep=0cm, column sep=2cm, ampersand replacement=\&]
{
  |(left)| $J_K(y,1+h)$
  \&
  |(mid)| $
    \begin{aligned}
      \frac{1}{\Delta_K(y)}
      +\sum_{j=1}^\infty
      \frac{P^{(j)}_K(y)}{\Delta_K^{2j+1}(y)}\,h^j
    \end{aligned}
  $
  \&
  |(right)| $
    \begin{aligned}
      \frac{F_K(x,1+h)}{x^{1/2}-x^{-1/2}}
    \end{aligned}
  $
  \\
};

\draw[<-, shorten <=2pt, shorten >=2pt]
  (left.east) -- node[above] {$y \to 0$}
                node[below] {$x \to 1$}
  (mid.west);

\draw[->, shorten <=2pt, shorten >=2pt]
  (mid.east) -- node[above] {$y\to \infty$}
               node[below] {$x \to 0$}
  (right.west);

\end{tikzpicture}

\end{center}

There is an important difference between computing $P_K^{(j)}$ from $J_K$ versus from $F_K$: The Habiro series reads $P_K^{(j)}$ `from the left'---it first determines the constant coefficient in $y$, then the linear, etc.---whereas $F_K$ reads the coefficients of $P_K^{(j)}$ from the right, starting with the leading coefficient. This behavior can be illustrated as follows:

\begin{center}
\begin{tikzpicture}
  \node at (0,0) {$P_K^{(j)}(y) = a_0 y^0 + a_1 y^1 + \dots + a_{n-1} y^{n-1} + a_n y^n$};
  
  \draw[->] (-2.2,-.5) -- (-0.5,-.5) node[midway, below] {$J_K$};
  
  \draw[->] (3.9,-.5) -- (2.2,-.5) node[midway, below] {$F_K$};
\end{tikzpicture}
\end{center}

As a result, we obtain the following corollary of \cref{thm:L-P1}:

\begin{corollary}\label{cor:L-alex-habiro} Let $K$ be a nice knot. Then,
\begin{equation}\label{eq:L-alex-habiro}
\el(K) = \sum_{i=0}^{2\dK-1} \frac{1}{i!}\frac{\partial^i\Delta_K^3(y)}{\partial y^i}\biggr|_{y=0} \frac{\partial \an{2\dK-i}{K}}{\partial q}\biggr|_{q=1}.
\end{equation}
\end{corollary}

\begin{proof}

Since $\prod_{i=1}^j\left(y+2-(1+h)^i-(1+h)^{-i}\right)=y^j+O(h^2)$, the Habiro series evaluated at $q=1+h$ and expanded at $h=0$ takes the form
\begin{equation*}
    J_K(y,1+h)=\sum_{j=0}^\infty y^j\left(\anq{j}{K}{1}+h\frac{\partial \an{j}{K}}{\partial q}\biggr|_{q=1}\right) + O(h^2).
\end{equation*}
From the relation between the Habiro series and the MMR expansion we deduce that
\begin{equation*}
    \sum_{j=0}^\infty y^j\frac{\partial \an{j}{K}}{\partial q}\biggr|_{q=1}= \frac{P^{(1)}_K(y)}{\Delta_K^3(y)}.
\end{equation*}
By \cref{thm:L-P1}, the coefficient of $y^{2\dK}$ in $P^{(1)}_K(y)$ is $\el(K)$. Thus the coefficient of $y^{2\dK}$ in
\begin{equation*}
    \sum_{j=0}^\infty \Delta_K^3(y)\, \frac{\partial \an{j}{K}}{\partial q}\biggr|_{q=1}\, y^j
\end{equation*}
is $\el(K)$, which gives us the desired result.
\end{proof}

\cref{cor:L-alex-habiro} together with \cref{thm:L-Hopf-homogeneous} imply an explicit formula for Hopf invariant in terms of the Habiro coefficients $a_n(q), n=0,1,2,\dots$ and the Alexander polynomial $\Delta_K(x)$: 

\begin{theorem}\label{thm:MMR-1-homogeneous} Let $K$ be a homogeneous braid knot, or a fibered knot of at most 12 crossings. Then,
    \begin{equation}\label{eq:hopf}
        \hopfK = \genusK-\sum_{i=0}^{2\genusK-1} \frac{1}{i!}\frac{\partial^i\Delta_K^3(y)}{\partial y^i}\biggr|_{y=0} \frac{\partial \an{2\genusK-i}{K}}{\partial q}\biggr|_{q=1}.
    \end{equation}
\end{theorem}
Note that \cref{conj:L-hopf} would imply the formula \eqref{eq:hopf} for all fibered knots.

\subsection{General roots of unity}\label{ss:mmr_roots_of_unity}

Let $\mathcal{Z} \subset \C$ denote the set of all rational roots of unity. Habiro conjectured a version of the MMR expansion at any root of unity:
\begin{conjecture}\cite[Conj. 7.4]{Habiro2007unified}\label{conj:Habiro-ADO-MMR} Let $\omega$ be a primitive $p$-th root of unity for $p \geq 1$ and let $y=x+x^{-1}-2$. Let 
\[
    \gamma_{\omega}:\Lambda\rightarrow \Z[\omega][[q-\omega,y]]
\]  
be the natural homomorphism induced by $\Z[q,q^{-1},y]\subset \Z[\omega][q,q^{-1},y]$. Then,
\[
    \gamma_{\omega}(J_K)= 
    \sum_{j \geq 0} \frac{P_{K,\omega}^{(j)}(x)}{\Delta_K(x^p)^{2j+1}}(q-\omega)^j
\]
for some polynomials $P_{K,\omega}^{(j)}(x) \in \Z[\omega][x+x^{-1}]$. 
\end{conjecture}

We expect a similar structure for the $F_K$ series of nice knots.

\begin{conjecture}\label{conj:strong-FK-ADO} Let $\omega$ be a primitive $p$-th root of unity for $p \geq 1$ and let $K$ be a nice knot. Then the series $F_K(x,q)$ admits an expansion at $q=\omega$ of the form
\begin{equation}\label{eq:MMR_roots_FK}
    \Fu_K(x,h+\omega) = \sum_{j\geq 0}\frac{P_{K,\omega}^{(j)}(x)}{\Delta_K^{2j+1}(x^p)} h^j,
\end{equation}
where $P_{K,\omega}^{(j)}(x)\in\Z[\omega][x+x^{-1}]$ coincide with the polynomials in \cref{conj:Habiro-ADO-MMR}. Both sides of \eqref{eq:MMR_roots_FK} are understood as power series in $x$ and $h$.
\end{conjecture}

We write $P_{K,p}^{(j)}(x)$ for $P_{K,\zeta_p}^{(j)}(x)$, where $\zeta_p=e^{2\pi i/p}$. 

\begin{remark}
    Habiro's conjecture \ref{conj:Habiro-ADO-MMR} is valid at order $0$ in $h=q-\omega$, by work of Willetts and Beliakova--Hikami \cite{Willetts2022,Beliakova2021}. More specifically, they proved that $P_{K,p}^{(0)}(x)=N_p(K,x)$, where $N_p(K,x)$ are the ADO polynomials. Hence \cref{conj:strong-FK-ADO} generalizes \cref{conj:FK-ADO} to higher orders in $h$.
\end{remark}

Conjecture \ref{conj:strong-FK-ADO} would have consequences analogous to Theorem \ref{thm:L-P1}, and Propositions \ref{prop:deg-Pj-MMR-posL} and \ref{prop:deg-Pj-MMR-negL} for the polynomials $P_{K,\omega}^{(j)}(x)$, generalizing the case $j=0$ in  \cite[Thm.\!1]{NeumannVanDerVeen24}:

\begin{proposition}\label{conj:deg-Pj-MMR-ADO}
    Assuming \cref{conj:strong-FK-ADO} for a nice knot $K$, we have
\begin{enumerate}
    \item If $\el(K)\geq 0$, then 
    \(\deg_x P^{(j)}_{K,\omega}(x)=-((2j+1)p-1)\dK\)
    with leading coefficient $\binom{\el(K)}{j}$ $\omega^{\el(K)-j}$ for $0\leq j\leq g(K)$, and 
    \(\deg_x P^{(j)}_{K,\omega}(x)>-((2j+1)p-1)\dK\)
    for $j>\el(K)$.
    \item If $\el(K)<0$, then
    \(\deg_x P^{(j)}_{K,\omega}(x)=-((2j+1)p-1)\dK\)
    with leading coefficient $(-1)^j$ $\binom{-\el(K)+j-1}{j}\omega^{\el(K)-j}$ for all $j\geq 0$.
\end{enumerate}
In particular, $\el(K)$ is the modulus of the coefficient of $x^{2(3p-1)\dK}$ in $P^{(1)}_{K,\omega}(x)$.
\end{proposition}

\begin{proof}
We take the expansion of the $F_K$ invariant near $\omega$, modulo overall sign, and consider the leading term:
\begin{equation*}
\begin{split}
    \lt(\Fu_K(x,h+\omega))&=(h+\omega)^{\el(K)} x^{\dK}=
    \begin{cases}
         x^{\dK} \sum_{j=0}^{\el(K)} \binom{\el(K)}{j} \omega^{\el(K)-j}h^j & \text{ if } \el(K) \geq 0,\\
         x^{\dK} \sum_{j=0}^{\infty} (-1)^j \binom{-\el(K)+j-1}{j} \omega^{\el(K)-j}h^j & \text{ if } \el(K) <0.
    \end{cases}
\end{split}
\end{equation*}
The claim follows from comparison with \eqref{eq:MMR_roots_FK}, as in the proofs of \cref{prop:deg-Pj-MMR-posL,prop:deg-Pj-MMR-negL}.
\end{proof}

We end this section with two examples of $P^{(1)}_{K,p}(x)$ for the right-handed trefoil and figure-eight knot. For examples of $P_{K,p}^{(0)}(x)$, see \cite{GukovNakajima2021, Chae2020}. Note that we only write terms with nonnegative $x$-powers. We verified that these polynomials appear in the expansions of both $F_K$ and $J_K$. 

For trefoil, we have $\el(3_1^r)=1$ and $g(3_1^r)=1$, so the $x$-degree of $P_{3_1^r,p}^{(1)}$ is $-(3p+1)$, with the leading coefficient $\el(3_1^r)\zeta_p^{\el(3_1^r)-1} = 1$. The polynomials $P_{3_1^r,p}^{(1)}$ read:

\begin{center}
\begin{tabular}{c  c}
    $p$ & $P^{(1)}_{3_1^r,p}(x)$ \text{with} $\zeta=\zeta_p$ \text{(nonnegative part)}\\
    \midrule[0.1pt]
    3 &$(-3 + 6 \zeta) - 8 x - 8 x^2 + (-5 + \zeta) x^3 +(1+\zeta) x^4 +(1+\zeta) x^5 + (1 - 2 \zeta) x^6 + x^7 + x^8$ \\
    \midrule[0.1pt]
    4 & \makecell[c]{$(-20 + 12 \zeta) - (15 - 12 \zeta) x - 15 x^2 - (15 - 2 \zeta) x^3 - (8 - 2 \zeta) x^4 + (1 + 2\zeta) x^5+$\\$ (1 + 4 \zeta) x^6 + (1 - 2 \zeta) x^7 + (4 - 2 \zeta) x^8 + (1 - 2 \zeta) x^9 + x^{10} + x^{11} $}\\
    \midrule[0.1pt]
    5 & \makecell[c]{$(-12 + 19 \zeta + 21 \zeta^2) + (-24 + 19 \zeta + 21 \zeta^2) x + (-24 + 
    19 \zeta) x^2 + (-24 + 4 \zeta) x^3 +$\\ $(-24 + 
    4 \zeta + \zeta^2) x^4 + (-12 + 4 \zeta + \zeta^2) x^5 + (4 + 4 \zeta + \zeta^2) x^6 + (4 + 4 \zeta + 5 \zeta^2) x^7+$\\ $(4 - 2 \zeta + 5 \zeta^2) x^8 + (4 - 2 \zeta - 3 \zeta^2) x^9 + (1- 2 \zeta - 3 \zeta^2) x^{10} + (1 - 2 \zeta - 
    3 \zeta^2) x^{11} + $ \\ $(1 - 2 \zeta) x^{12} + x^{13} + x^{14}$}
\end{tabular}
\end{center}

The figure-eight knot $4_1$ is amphichiral, hence $\el(4_1)=0$. As $\el(4_1)=0$ and $g(4_1)=1$, the degree of the polynomial $P_{4_1,p}^{(1)}$ is strictly larger than $-3p+1$, since $\el(4_1)\zeta_p^{\el(4_1)-1} = 0$.  The polynomials $P_{4_1,p}^{(1)}$ read:

\begin{center}
\begin{tabular}{c c}
    $p$ & $P^{(1)}_{4_1,p}(x)$ \text{with} $\zeta=\zeta_p$ \text{(nonnegative part)}\\
    \midrule[0.1pt]
    3 & $(8 - 8 \zeta) + (-1 + \zeta) x + (3 - 3 \zeta) x^2 + (-5 + 
    5\zeta) x^3 + (-1 + \zeta) x^5 + (1 - \zeta) x^6$\\
    \midrule[0.1pt]
    4 &$ 54 + 18 x + 8 x^2 - 32 x^4 - 12 x^5 - 4 x^6 - 2 x^7 + 6 x^8 + 2 x^9$ \\
    \midrule[0.1pt]
    5 & \makecell[c]{$(44 + 33 \zeta - 33 \zeta^2 - 44 \zeta^3) + (35 + 33 \zeta - 33 \zeta^2 - 35 \zeta^3)x + (17 - 3 \zeta + 3 \zeta^2 - 17 \zeta^3)x^2 +$ \\
   $ (17 \zeta - 17 \zeta^2)x^3  + (-16 - 15 \zeta + 15 \zeta^2 + 16 \zeta^3) x^4 
    + (-24 - 18 \zeta + 18 \zeta^2 + 24 \zeta^3)x^5 +$\\
    $ (-19 - 21 \zeta + 21 \zeta^2 + 19 \zeta^3)x^6 
    +  (-9 - \zeta + \zeta^2 + 9 \zeta^3)x^7 + (-2 - 6 \zeta + 6 \zeta^2 + 2 \zeta^3)x^8 +$ \\
    $(2 + 2 \zeta - 2 \zeta^2 - 2 \zeta^3)x^9 + (4 + 3 \zeta - 3 \zeta^2 - 4 \zeta^3)x^{10} + (3 + 4 \zeta - 4 \zeta^2 - 3 \zeta^3) x^{11} + $\\
    $(1 - \zeta^3) x^{12}$}
\end{tabular}
\end{center}

\section{Slopes}\label{sec:slopes}

The $F_K$ slopes of fibered Montesinos knots up to 10 crossings are collected in the tables below. Here we give a detailed example of a slope computation for $K=12n_{242}$, also known as the $(-2,3,7)$-pretzel knot. The $F_K$ reads
    \begin{equation*}
    \begin{split}
    F_K = &-q^5 x^4+q^7 x^7+q^9 x^9-q^{10} x^{10}+q^{11} x^{11}-q^{13} x^{12}+q^{14} x^{13}-(q^{15}+q^{16}) x^{14}\\
    &+(q^{17}+q^{18}) x^{15}-2 q^{19} x^{16}+(q^{20}+q^{21}+q^{22})
   x^{17} -(q^{22}+q^{23}+q^{24}) x^{18}\\
   &+(2 q^{25}+q^{26}) x^{19}-(q^{26}+2 q^{27}+q^{28}+q^{29}) x^{20} +(q^{28}+q^{29}+2 q^{30}+q^{31}) x^{21} - \cdots\\
   \end{split}
    \end{equation*}

The sequence of the maximal $q$-degrees (starting at $x^{10}$), and its differences, are
\begin{equation*}
\begin{split}
\delta_{\text{max}}:& 10, 11, 13, 14, 16, 18, 19, 22, 24, 26, 29, 31, 34, 37, 39, 43, 46, 
49, 53, 56, 60, 64, 67, 72, \dots\\
\Delta \delta_{\text{max}}:& 1, 2, 1, 2, 2, 1, 3, 2, 2, 3, 2, 3, 3, 2, 4, 3, 3, 4, 3, 4, 4, 3, 5, \dots\\
\Delta^2 \delta_{\text{max}}:& 1, -1, 1, 0, -1, 2, -1, 0, 1, -1, 1, 0, -1, 2, -1, 0, 1, -1, 1, 0, -1, 2, \dots
\end{split}
\end{equation*}
We observe that $\Delta^2 \delta_{\text{max}}$ is periodic. We take its generating function, and using \cite[Lemma 4.4]{garoufalidis2010slopes}, we obtain the generating function for the original sequence $\delta_{\text{max}}$:
\[
\frac{t \frac{t \frac{-1 + t - t^2 + t^4 - 2 t^5 + t^6}{-1 + t^8} + 1}{1 - t} + 10}{1 - t} = \frac{-10 + 9 t - t^2 + t^3 - t^4 + t^6 - 2 t^7 + 11 t^8 - 
 9 t^9}{(-1 + t)^3 (1 + t + t^2 + t^3 + t^4 + t^5 + t^6 + t^7)}.
\]
From the generating function, we deduce that the slope is $16$ using \cite[Cor\! 4.3]{garoufalidis2010slopes}.

Similarly, the sequence of minimal $q$-degrees is
\[
\delta_{\text{min}}: 10, 11, 13, 14, 15, 17, 19, 20, 22, 25, 26, 28, 31, 33, 35, 38, 41, 43, 46, 50, 52, 55, 59, 62, 65, \dots
\]
and we find their generating function
\[
\frac{-11 + 9 t + t^2 - t^4 + t^6 - t^7 - t^8 + 2 t^9 + 10 t^{10} - 
 10 t^{11}}{(-1 + t)^3 (1 + t + t^2 + t^3 + t^4 + t^5 + t^6 + t^7 + 
   t^8 + t^9)},
\]
from which we read off the slope $20$. 

The $(-2,3,7)$-pretzel knot admits many exceptional surgeries, namely $16,17,37/2,18,19$ and $20$. The $F_K$ slopes are exactly the minimal and maximal exceptional slope, while the Jones slopes are $0$ and $37/2$. The slopes $16,37/2$ and $20$ are all toroidal. Finally note that the width of the $f_n$ polynomial is approximately $(1/16-1/20)n^2$, hence $f_n$ is roughly 370 times shorter than the colored Jones polynomial $J_n$.

\begin{remark}
Positive Montesinos knots have nonnegative boundary slopes \cite{Ichihara2008}.
The knot $10_{145}$ is strongly quasi-positive, but not positive. Interestingly, $F_K$ detects the negative slope $-2$, witnessing non-positivity, while the Jones slope is $0$.
\end{remark}

\begin{remark}
One may ask what are the slopes of the inverted Habiro coefficients $a_{-n}(q)$ defined by Park, which conjecturally extend the usual Habiro coefficients $a_n(q)$ to negative integers. Note that the coefficients $a_{-n}(q)$, $n=1,2,3,\dots$ are equivalent to the coefficients $f_j(q), j=0,1,2,\dots$ \cite{Park21}. The coefficients $a_{-n}(q)$ should correspond to a slope of the Newton polygon of the $C$-polynomial, defined in \cite{GaroufalidisSun2005}. Interestingly, for most fibered knots we considered, the leading term of the growth of minimal powers of $a_{-n}(q)$ is $-1/2$. This is in particular true for all knots in KnotInfo whose braid representative is homogeneous (up to 13 crossings) except the figure-eight knot $4_1$, for which $a_{-n}=1$ for all $n$.
\end{remark}

We give a table of slopes of the $F_K$ invariant of fibered Montesinos knots up to 10 crossings\footnote{For Montesinos knots, all boundary slopes are known. Otherwise this condition plays no role for our considerations.}, and include the Jones slopes for comparison. We write the $F_K$ slopes in bold, and the Jones slopes are underlined.

The exceptional slopes of a knot $K$ tend to span an interval bounded by two boundary slopes. See \cite[Thm.\! 1]{Ichihara2023} for a precise statement in the case of Montesinos and alternating knots. We include the boundaries of this interval, denoted by an asterisk, for comparison with the $F_K$ and Jones slopes.

\renewcommand{\arraystretch}{1.1}

\begin{longtable}{llll}
Name & Boundary Slopes  & Name & Boundary slopes\\
\hline
$3_{1} = \mathrm{K}_{1/3}$ &   $\jsc{0}$, $\jsc{\fsc{6}}$ & 
$4_{1} = \mathrm{K}_{2/5}$ &   $\exc{\jsc{\fsc{-4}}}$, $0$, $\exc{\jsc{\fsc{4}}}$\\ 
$5_{1} = \mathrm{K}_{1/5}$ &   $\jsc{0}$, $\jsc{\fsc{10}}$ &
$6_{2} = \mathrm{K}_{4/11}$ &   $\jsc{\fsc{-4}}$, $0$, $2$, $\exc{\jsc{\fsc{8}}}$ \\ 
$6_{3} = \mathrm{K}_{5/13}$ &   $\jsc{\fsc{-6}}$, $-2$, $0$, $2$, $\jsc{\fsc{6}}$ &
$7_{1} = \mathrm{K}_{1/7}$ &   $\jsc{0}$, $\jsc{\fsc{14}}$ \\
$7_{6} = \mathrm{K}_{7/19}$ &   $\jsc{\fsc{-4}}$, $0$, $\fsc{4}$, $6$, $\jsc{10}$ &
$7_{7} = \mathrm{K}_{8/21}$ &   $\jsc{-8}$, $\fsc{-4}$, $0$, $\jsc{\fsc{6}}$ \\
$8_{2} = \mathrm{K}_{6/17}$ &   $\jsc{\fsc{-4}} $, $0$, $6$, $\exc{\jsc{\fsc{12}}}$ & 
$8_{7} = \mathrm{K}_{9/23}$ &   $\jsc{\fsc{-10}}$, $-6$, $-2$, $0$, $\jsc{\fsc{6}}$ \\
$8_{9} = \mathrm{K}_{7/25}$ &   $\jsc{\fsc{-8}}$, $-2$, $0$, $2$, $\jsc{\fsc{8}}$ &
$8_{12} = \mathrm{K}_{12/29}$ &   $\jsc{-8}$, $\fsc{-4}$, $0$, $\fsc{4}$, $\jsc{8}$ \\
$9_{1} = \mathrm{K}_{1/9}$ &   $\jsc{0}$, $\jsc{\fsc{18}}$ &
$9_{11} = \mathrm{K}_{14/33}$ &   $\jsc{-14}$, $-10$, $\fsc{-4}$, $0$, $\jsc{\fsc{4}}$ \\
$9_{17} = \mathrm{K}_{14/39}$ &   $\jsc{-8}$, $\fsc{-4}$, $-2$, $0$, $4$, $\jsc{\fsc{10}}$ &
$9_{20} = \mathrm{K}_{15/41}$ &   $\jsc{\fsc{-4}}$, $0$, $2$, $\fsc{6}$, $8$, $\jsc{14}$ \\ 
$9_{26} = \mathrm{K}_{18/47}$ &  $\jsc{-12}$, $-8$, $-6$, $\fsc{-4}$, $0$, $\jsc{\fsc{6}}$ & 
$9_{27} = \mathrm{K}_{19/49}$ &   $\jsc{\fsc{-8}}$, $-4$, $-2$, $0$, $2$, $\fsc{4}$, $6$, $\jsc{10}$ \\
$9_{31} = \mathrm{K}_{21/55}$ &   $\jsc{\fsc{-6}}$, $-2$, $0$, $2$, $\fsc{6}$, $\jsc{12}$ & 
$10_{2} = \mathrm{K}_{8/23}$ &   $\jsc{\fsc{-4}}$, $0$, $10$, $\exc{\jsc{\fsc{16}}}$ \\ 
$10_{5} = \mathrm{K}_{13/33}$ &   $\jsc{\fsc{-14}}$, $-10$, $-6$, $0$, $\jsc{\fsc{6}}$ &
$10_{9} = \mathrm{K}_{11/39}$ &   $\jsc{\fsc{-8}}$, $0$, $2$, $6$, $\jsc{\fsc{12}}$\\ 
$10_{17} = \mathrm{K}_{9/41}$ &   $\jsc{\fsc{-10}}$, $-2$, $0$, $2$, $\jsc{\fsc{10}}$ &
$10_{29} = \mathrm{K}_{26/63}$ &   $\jsc{-8}$, $\fsc{-4}$, $-2$, $0$, $2$, $\fsc{4}$, $8$, $\jsc{12}$\\ 
$10_{41} = \mathrm{K}_{26/71}$ &   $\jsc{-8}$, $\fsc{-4}$, $0$, $2$, $\fsc{4}$, $6$, $8$, $\jsc{12}$ & 
$10_{42} = \mathrm{K}_{31/81}$ &   $\jsc{-10}$, $-6$, $\fsc{-4}$, $-2$, $0$, $2$, $\fsc{4}$, $6$, $\jsc{10}$\\ 
$10_{43} = \mathrm{K}_{27/73}$ &   $\jsc{-10}$, $-6$, $\fsc{-4}$, $0$, $\fsc{4}$, $6$, $\jsc{10}$ &
$10_{44} = \mathrm{K}_{30/79}$ &   $\jsc{-8}$, $\fsc{-4}$, $0$, $2$, $4$, $\fsc{6}$, $\jsc{12}$\\ 
$10_{45} = \mathrm{K}_{34/89}$ &   $\jsc{-10}$, $-6$, $\fsc{-4}$, $-2$, $0$, $2$, $\fsc{4}$, $6$, $\jsc{10}$ &
\end{longtable}

\begin{longtable}{ll}
Name & Boundary Slopes \\
\hline
$8_{10} = \mathrm{K}(1/3, 2/3, 1/2)$ & $\jsc{-10}$, $-8$, $\fsc{-6}$, $0$, $2$, $\jsc{\fsc{6}}$\\ 
$8_{19} = \mathrm{K}(1/3, 1/3, -1/2)$ &  $\jsc{0}$, $\jsc{\fsc{12}}$\\ 
$8_{20} = \mathrm{K}(1/3, 2/3, -1/2)$ &  $\jsc{\fsc{-10}}$, $0$, $\jsc{\fsc{8/3}}$ \\ 
$8_{21} = \mathrm{K}(2/3, 2/3, -1/2)$ &  $\jsc{\fsc{-1}}$, $0$, $2$, $6$, $ \jsc{\fsc{12}}$\\ 
$9_{22} = \mathrm{K}(3/5, 1/3, 1/2)$ &   $\jsc{-8}$, $\fsc{-4}$, $-2$, $0$, $2$, $4$, $\fsc{6}$, $8$, $\jsc{10}$\\ 
$9_{24} = \mathrm{K}(1/3, 2/3, 3/2)$ &   $\jsc{-8}$, $-6$, $\fsc{-4}$, $-2$, $0$, $\fsc{4}$, $6$, $\jsc{10}$\\ 
$9_{28} = \mathrm{K}(2/3, 2/3, 3/2)$ &   $\jsc{-6}$, $\fsc{-4}$, $-2$, $0$, $2$, $\fsc{6}$, $8$, $\jsc{12}$\\ 
$9_{30} = \mathrm{K}(3/5, 2/3, 1/2)$ &   $\jsc{-10}$, $-6$, $\fsc{-4}$, $-2$, $0$, $2$, $4$, $\fsc{6}$, $\jsc{8}$\\ 
$9_{36} = \mathrm{K}(2/5, 1/3, 1/2)$ &   $\jsc{-14}$, $-12$, $-10$, $-8$, $-6$, $\fsc{-4}$, $-2$, $0$, $\jsc{\fsc{4}}$\\ $9_{42} = \mathrm{K}(2/5, 1/3, -1/2)$ &   $\exc{\jsc{-6}}$, $\exc{\fsc{-8/3}}$, $0$, $\jsc{\fsc{8}}$\\ $9_{43} = \mathrm{K}(3/5, 1/3, -1/2)$ &   $\jsc{\fsc{-4}}$, $0$, $6$, $\fsc{8}$, $\jsc{32/3}$\\ 
$9_{44} = \mathrm{K}(2/5, 2/3, -1/2)$ &   $\jsc{-14/3}$, $-2$, $\fsc{-1}$, $0$, $2$, $\jsc{\fsc{10}}$\\
$9_{45} = \mathrm{K}(3/5, 2/3, -1/2)$ &   $\jsc{-14}$, $-10$, $-8$, $\fsc{-4}$, $-2$, $0$, $\jsc{\fsc{1}}$\\
$9_{48} = \mathrm{K}(2/3, 2/3, -1/3)$ &   $\jsc{-11}$, $-8$, $\fsc{-4}$, $0$, $\jsc{\fsc{4}}$\\ 
$10_{46} = \mathrm{K}(1/5, 1/3, 1/2)$ &   $\jsc{\fsc{-4}}$, $0$, $2$, $6$, $\fsc{8}$, $12$, $14$, $\jsc{16}$\\ 
$10_{47} = \mathrm{K}(1/5, 2/3, 1/2)$ &   $\jsc{-14}$, $-12$, $-10$, $\fsc{-6}$, $-4$, $0$, $2$, $\jsc{\fsc{6}}$\\ $10_{48} = \mathrm{K}(4/5, 1/3, 1/2)$ &   $\jsc{\fsc{-10}}$, $-6$, $-4$, $-2$, $0$, $2$, $\fsc{6}$, $8$, $\jsc{10}$ \\
$10_{59} = \mathrm{K}(2/5, 3/5, 1/2)$ &   $\jsc{-8}$, $\fsc{-4}$, $0$, $2$, $\fsc{4}$, $6$, $8$, $10$, $\jsc{12}$\\ 
$10_{60} = \mathrm{K}(3/5, 3/5, 1/2)$ &   $\jsc{-12}$, $-8$, $\fsc{-4}$, $0$, $2$, $\fsc{6}$, $\jsc{8}$\\ 
 
$10_{62} = \mathrm{K}(1/4, 1/3, 2/3)$ &   $\jsc{-14}$, $-10$, $\fsc{-8}$, $-2$, $0$, $\jsc{\fsc{6}}$\\ 
$10_{64} = \mathrm{K}(3/4, 1/3, 1/3)$ &   $\jsc{\fsc{-8}}$, $-2$, $0$, $2$, $4$, $6$, $\fsc{8}$, $\jsc{12}$\\ 
$10_{69} = \mathrm{K}(3/5, 2/3, 2/3)$ &  $\jsc{-14}$, $-12$, $-10$, $-8$, $-6$, $\fsc{-4}$, $-2$, $0$, $2$, $\jsc{\fsc{6}}$\\ 
$10_{70} = \mathrm{K}(2/5, 1/3, 3/2)$ &   $\jsc{-12}$, $-10$, $-8$, $-6$, $\fsc{-4}$, $-2$, $0$, $2$, $\fsc{4}$, $\jsc{8}$\\ $10_{71} = \mathrm{K}(2/5, 2/3, 3/2)$ &   $\jsc{-10}$, $-6$, $\fsc{-4}$, $-2$, $0$, $2$, $\fsc{4}$, $6$, $8$, $\jsc{10}$\\ 
$10_{73} = \mathrm{K}(3/5, 2/3, 3/2)$ &   $\jsc{-14}$, $-10$, $-8$, $-6$, $\fsc{-4}$, $-2$, $0$, $2$, $\fsc{4}$, $\jsc{6}$\\ 
$10_{75} = \mathrm{K}(2/3, 2/3, 5/3)$ &  $\jsc{\fsc{-8}}$, $-2$, $0$, $\fsc{4}$, $8$, $\jsc{12}$\\ $10_{78} = \mathrm{K}(2/3, 2/3, 5/2)$ &   $\jsc{\fsc{-4}}$, $-2$, $0$, $\fsc{4}$, $6$, $10$, $12$, $\jsc{16}$\\ 
$10_{124} = \mathrm{K}(1/5, 1/3, -1/2)$ &   $\jsc{0}$, $\jsc{\fsc{15}}$\\ 
$10_{125} = \mathrm{K}(1/5, 2/3, -1/2)$ &   $\fsc{-32/5}$, $-4$, $0$, $\jsc{\fsc{10}}$\\ $10_{126} = \mathrm{K}(4/5, 1/3, -1/2)$ &  $\jsc{-14}$, $\fsc{-8}$, $-6$, $-4$, $0$, $\fsc{8/3}$\\ $10_{127} = \mathrm{K}(4/5, 2/3, -1/2)$ &   $\fsc{-1}$, $\jsc{0}$, $2$, $6$, $\fsc{8}$, $10$, $\jsc{16}$\\ $10_{132} = \mathrm{K}(2/7, 1/3, -1/2)$ &   $\jsc{\fsc{-14}}$, $\exc{-2}$, $0$, $\exc{\fsc{3/2}}$\\ 
$10_{133} = \mathrm{K}(2/7, 2/3, -1/2)$ &   $\fsc{-1/2}$, $\jsc{0}$, $2$, $10/3$, $6$, $\jsc{\fsc{16}}$\\ 
$10_{139} = \mathrm{K}(1/4, 1/3, -2/3)$ &   $\jsc{0}$, $\fsc{\exc{12}}$, $\exc{13} $, $\jsc{18}$, $\fsc{20}$\\ 
$10_{140} = \mathrm{K}(1/4, 1/3, -1/3)$ &  $\exc{\jsc{\fsc{-8/5}}}$, $\exc{0}$, $\fsc{14}$\\ 
$10_{141} = \mathrm{K}(1/4, 2/3, -1/3)$ &   $\fsc{-9/2}$, $-2$, $0$, $2$, $4$, $\jsc{\fsc{12}}$\\
$10_{143} = \mathrm{K}(3/4, 1/3, -1/3)$ &  $\jsc{-14}$, $-8$, $\fsc{-6}$, $-2$, $0$, $\fsc{8/3}$\\ 
$10_{145} = \mathrm{K}(2/5, 1/3, -2/3)$ &  $\jsc{-18}$, $\exc{-6}$, $\exc{\fsc{-4}}$, $\jsc{0}$, $\fsc{2}$\\
\end{longtable}

\section{Beyond fibered knots}\label{sec:nonfibered}

The last section is experimental in nature, as there is no available definition of the $F_K$ invariant of non-fibered knots. So far, there were only isolated examples of the $F_K$ series for non-fibered knots, constructed by Park using an extension of the state sum model (described below), or by interpolating explicit formulas of the colored Jones polynomials of twist knots from \cite{Lovejoy2021}. 

We extend the pool of examples and collect new observations. Among them is the relation between the $F_K$ invariant and the stability series of the colored Jones polynomials, and the failure of the previously anticipated relation between the $F_K$ invariant and the MMR expansion. The most striking observation, however, is the similarity between the behavior of the $F_K$ invariant and that of knot Floer homology, which generalizes our previous discussion relating $F_K$ to the Hopf invariant beyond the fibered case.

The main new feature in the non-fibered case is that the coefficients of the $F_K$ series are expected to be infinite $q$-series. This was originally motivated by the relation between $F_K$ and the Alexander polynomial: For non-fibered knots, the Alexander polynomial is not necessarily monic, and hence the expansion of the inverse Alexander polynomial may have rational coefficients. If we expect \eqref{eq:FK-MMR} to hold, given that the coefficients of $F_K(x,q)$ must be integers, this rules out the possibility that its $x$-coefficients are finite polynomials in $q$. Nevertheless, we find that even in cases where the Alexander polynomial is monic, or where the limit of the $F_K$ invariant at $q=1$ is divergent (and the relation \eqref{eq:FK-MMR} fails), the $F_K$ invariant of a non-fibered knot is still an infinite Laurent power series in $q$, which strengthens our expectation.

\subsection{Strongly quasipositive knots}\label{sec:sqp}

\begin{figure}
\includegraphics[]{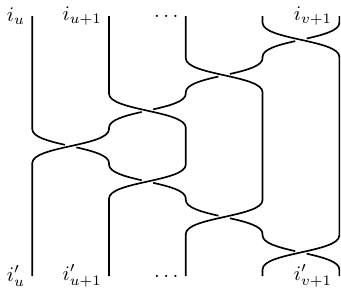}
\caption{The braid $\sigma_{u,v}$ and its labels.}
\label{fig:band-generator}
\end{figure}

We recall the notion of strongly quasipositive knots. 
For $u, v \in \{1,\dots, n-1\}$ with $u\leq v$, we consider the element
\[
    \sigma_{u,v} = (\sigma_v \dots \sigma_{u+1}) \sigma_{u} (\sigma_{v} \dots \sigma_{u+1})^{-1} \in B_n,
\]
commonly known as \emph{band generator}.
It corresponds to a braid that switches the $u$ and $v$ positions positively (see \cref{fig:band-generator}). In our convention, $\sigma_{u,u}$ is the usual Artin generator $\sigma_u$. A knot is \emph{strongly quasipositive} if it can be written as the closure of a braid $\beta$ composed from band generators:
\[
    \beta = \prod_k \sigma_{u_k,v_k}.
\]

In this section, we consider mirror images of strongly quasipositive knots, which we call \emph{strongly quasinegative}. This is necessary due to our conventions on Jones polynomial. In particular, if one uses the opposite convention (as in \cite{GaroufalidisLe2015}), this section would apply to strongly quasipositive knots. 

\subsection{Stratified state sum}

In \cite{Park20}, Park constructed a candidate for the $F_K$ series of some non-fibered, strongly quasinegative knots, using the notion of the \emph{stratified} state sum, which we now review. Consider a braid $\beta$ representing the knot $K$, and equip $\beta$ with a constant inversion datum $\id\equiv -1$. Define the \emph{weight} $\eta(s) \in \N_{\geq 0}$ of a state $s \in \Omega(\id)$ as
\begin{equation}
    \eta(s)= -(n-1) + \sum_{j=2}^n |s(b_j)|,
\end{equation}
where the edges $b_1,\dots, b_n$ label the bottom segments of the braid (\cref{fig:braid_state_sum}). We shift by $-(n-1)$ so that the ground state $s_0 \equiv -1$ has weight $0$. Note that $\eta(s)> \eta(s_0)=0$ for any state $s\neq s_0$.

The weight function gives a decomposition (stratification) of the set of valid states $\Omega$:
\[
\Omega(\id) = \bigcup_{w=0}^{\infty} \Omega_w(\id)
\]
where $\Omega_{w}(\id)=\{s\in\Omega(\id) \mid \eta(s)=w\}$. The \emph{stratified state sum} of $\bd$ is defined as
\begin{equation}\label{eq:def-stratified-sum}
    Z^{\text{str}}(\bd)=\lim_{y\rightarrow 1}{\left(\sum_{w=0}^\infty  y^{-w}\sum_{s\in\Omega_w(\id)} P(s)\right)}.
\end{equation}
Whenever $Z^{\text{str}}(\bd)\in\FKuRing$, we define $F_K=(x^\frac12-x^{-\frac12})Z^{\text{str}}(\bd)$. As in the previous section, the unnormalized version is defined as $\tilde{F}_K(x,q)=\sum_{j\geq 0} \tilde{f}_j(q)x^{j+\dK}:=Z^{\text{str}}(\bd)$.

\subsubsection{Minimal contribution}

For nice knots, the only contribution to the leading term of the inverted state sum is given by the ground state, as proven in \cref{prop:ground-state}. In that case, by construction, the coefficients of $F_K$ in the variable $x$ are Laurent $q$-polynomials.

In contrast, for non-nice knots, we expect the leading term of the stratified state sum to be an infinite $q$-series, with infinitely many contributing states. However, we still expect the ground state $s_0$ to have a special contribution, which is studied in the following proposition.
\begin{proposition}\label{prop:stratified-ground-state} Let $K$ be a knot with strongly quasinegative braid representative $\bd\in B_n$. Then,
\begin{equation}
    P(s_0)=q^{-\genusK}x^{\genusK}.
\end{equation}
\end{proposition}

Define $\el(K)$ as the minimal $q$-exponent of the leading term of $F_K$:
\begin{equation*}      
\lt(F_K(x,q))=(\pm q^{\el(K)}+O(q^{\el(K)+1}))x^{\dK-1/2}.
\end{equation*}
Then, motivated by \cref{prop:stratified-ground-state}, we expect that $\dK=\genusK$ and $\el(K)=-\genusK$
for any strongly quasinegative knot for which the stratified state sum converges as a Laurent power series in $q$ and $x$. If $K$ is also fibered, then \cref{conj:L-hopf} would imply $\hopfK = -\genusK$, which can be proven using results in \cite{Rudolph1992,Rudolph1998,Giroux2006}.

\begin{proof}[Proof of \cref{prop:stratified-ground-state}] The ground state assigns $-1 $ to every edge. It follows from \cref{tab:xq-contributions} that 
\[
r_q^{\sgn(c),-1}=-r_x^{\sgn(c),-1}=\frac{\sgn(c)}{2}.
\]
Using $s(b_j)=-1$ for $j=2,\dots,n$, we obtain
\begin{equation}
    P(s_0) = \prod_{j=2}^n x^{-\frac12}q^{-\frac12-s(b_j)} \prod_{c\in V_\bd} q^{\frac{\sgn(c)}{2}}x^{-\frac{\sgn(c)}{2}} = q^{\frac{n-1}{2}+\frac{w(\beta)}{2}}x^{-\frac{n-1}{2}-\frac{w(\beta)}{2}}.
\end{equation} 

Let $\mir{\bd}$ be the mirror of $\bd$ given by changing the sign of every crossing. The closure of $\mir{\bd}$ is a strongly quasipositive knot. For any braid, the number of strands $n$, writhe $w$, and the genus of the braid closure are related by Bennequin's inequality. For quasipositive braids, the inequality is saturated \cite[Obs.\! 1.2.]{Ito2017}; i.e.,
   \begin{equation*}
       g(K)=-\frac{n-1}{2}+\frac{w(\mir{\bd})}{2}.
   \end{equation*}
   Since $w(\mir{\bd})=-w(\bd)$, the result follows.
\end{proof}

\subsection{\texorpdfstring{$F_K$}{FK} and the knot tail}

The colored Jones polynomial $\J{n}{K}$, normalized to start at 1, often stabilizes to a formal Laurent $q$-series. We call this series the \emph{tail} of the knot $K$, which was first studied in \cite{Dasbach_Lin_2006}. It has been shown to exist for alternating knots \cite{GaroufalidisLe2015} and, more generally, A-adequate knots\footnote{B-adequate in our convention.} \cite{Armond2013}.

We observed that many strongly quasinegative knots admit a tail. For all strongly quasinegative knots up to 10 crossings, the first non-zero $x$-coefficient $f_0(q)$ in \eqref{eq:F_K-terms}, computed from the stratified sum \eqref{eq:def-stratified-sum}, appears to converge to a Laurent power series in $q$ and recovers the tail $\Phi_{K,0}(q)$ via\footnote{The factor $(1-q)$ appears due to the different normalizations of the colored Jones polynomial and $F_K$, and disappears when considering the `unnormalized' analog of each invariant.}
\begin{equation}\label{eq:fk_tail}
    f_0(q) = -q^{-\genusK}(1-q)\Phi_{K,0}(q).
\end{equation}
In particular, we checked that this is true for all fibered strongly quasinegative knots up to 12 crossings: For every such knot, there exists a nice braid representative whose inversion datum assigns a negative integer to every segment. Therefore, \eqref{eq:fk_tail} implies that the tail of this subclass of knots is a monomial, up to the normalization factor $(1-q)^{-1}$.

\begin{remark} While it is known that negative knots, a subclass of strongly quasinegative knots, admit a tail (since they are B-adequate), this is not true for all strongly quasinegative knots: Indeed, the pretzel knot $P(5,-3,5)$ does not have a tail \cite{Lee22}, and it is strongly quasinegative by \cite{RudolphQuasipositivePretzels}\footnote{We thank Arpit Mittal for pointing out this fact to us.}.
\end{remark}

We gather the $f_0(q)$ of some strongly quasinegative knots in \cref{tab:tails}, computed using the stratified sum on their strongly quasipositive braid tabulated in \cite{knotinfo}; except for $8_{15}$, for which we used the braid $\sigma_1^2 \sigma_2^{-1}\sigma_1\sigma_3\sigma_2^3\sigma_3$, and $9_2$, for which we used the braid $\sigma_1\sigma_2\sigma_3\sigma_5\sigma_4\sigma_5\sigma_{1,4}$.

Given a knot with a tail, and once a suitable strongly quasinegative braid is found, the stratified sum is often a rather efficient way to compute the tail of strongly quasinegative knots (c.f. \eqref{eq:f0_infinite_intro}). It also gives closed formulas for the tail, some of which appear to be new. We gather some examples in \cref{tab:tails-closed-formulas}.
\renewcommand{\arraystretch}{1.1}
\begin{table}
    \centering
    \begin{tabular}{c l c}
        $\lc\Delta_K$ & $\pm q^{-\el(K)}f_{K,0}(q)$ & $K$ \\
        \hline
        $-1$ & $1 + 2 q + q^2 + q^3 + 2 q^4 + 2 q^5 + 2 q^6+\cdots$ & $\mathbf{12n_{801}}$ \\
        1 & $1$ & fibered, incl. $\mathbf{8_{19}}, \mathbf{10_{139,145,152,154,161}}$ \\
        1 & $1-q^2+2q^4-q^6-q^7+2q^9+\cdots$ & $\mathbf{13n_{1533}}$ \\
        2 & $1 - q + q^3 - q^6 + q^{10} - q^{15} + q^{21}-q^{28}+\cdots$ & $5_2, 7_3, 7_5, 9_{3,6,9,16}, \mathbf{10_{128,134,142}}$  \\
        2 & $1+q+q^4+q^9+q^{16}+q^{25}+q^{36}+\cdots$ & $\mathbf{12n_{293,321,750,830}}$\\
        3 & $1 - q + q^5 - q^8 + q^{16} - q^{21} + q^{33} + \dots$ & $7_2, 9_{4,7,47}, \mathbf{12n_{581}}$\\
        3 & $1 - 2 q + q^2 + 3 q^3 - 2 q^4+\cdots$ & $8_{15}, \mathbf{9_{49}}, 10_{49,66,80}$ \\
        4 & $1 - q + q^7 - q^{10} + q^{22} - q^{27} + q^{45} +\cdots$ & $9_2$ \\
        4 & $1 - 2 q + q^2 + 2 q^3-2q^4-q^6+\cdots$ & $7_4,9_{10,13,18,23}$ \\
        5 & $1 - 3 q + 3 q^2 + 4 q^3 - 8 q^4 - 2 q^5 + 2 q^6+\cdots$ & $9_{38}$ \\
        5 & $1 - 2 q + q^2 + q^3 - q^4 + 2 q^5 - 2 q^6+\cdots$ & $10_{55, 63}$ \\
        6 & $1 - 2 q + q^2 + q^3+\cdots$ & $9_5$ \\
        6 & $1 - 3 q + 3 q^2 + 3 q^3 - 7 q^4+\cdots$ & $10_{53}$\\
        7 & $1 - 2 q + q^2+3q^5-3q^6-q^7-q^8+\cdots$ &  $9_{35}$ \\
        7 & $1 - 3 q + 3 q^2 +2 q^3 -6 q^4 + 4 q^5+\cdots$ &  $10_{101}$\\
        8 & $1 - 4q + 6q^2 + 3q^3 - 17q^4 + \cdots$ & $10_{120}$ \\
    \end{tabular}
    \caption{Leading term of $F_K$, normalized to start at $1$, and leading coefficient of $\Delta_K$ for some strongly quasinegative knots, including all up to 10 crossings. Non-alternating knots are written in bold.} 
    \label{tab:tails}
\end{table}

\begin{remark}
 The \emph{stability series}\footnote{We use the letter $\Phi$ instead of original $F$ to avoid confusion with the GM series.} generalizes the tail of a knot $K$ \cite{GaroufalidisLe2015}. It is a Laurent power series in $q$ and $x$ of the form
\[
\Phi_K(x,q) = \Phi_{K,0}(q) + \Phi_{K,1}(q)\,x + \Phi_{K,2}(q)\,x^2 + \cdots
\]
where the first term $\Phi_{K,0}$ is the tail. For negative torus knots\footnote{Positive torus knots in the convention of \cite{GaroufalidisLe2015}.}, the stability series coincides with the $F_K$ series\footnote{
    If we use unnormalized versions $\Ju{n}{}$ and $\Fu_K(x,q)$, we have
    $\Fu_K(x,q)= (q^{-1} x)^\genusK \,\widetilde{\Phi}_K(x,q) \in \FKuRing$.
} \cite{GM}:
\begin{equation}\label{eq:stability-FK}
    F_K(x,q) = -(1-q) q^{-\genusK}x^{\genusK-1/2} \Phi_K(x,q),
\end{equation}
We observed that the stratified sum converges at every $x$-coefficient, and \eqref{eq:stability-FK} holds, for (suitable braids of) many strongly quasinegative knots, including $5_2,7_{3,5},9_{2,3,6,9,16,38,49}$ and $10_{49,63,128,134,142}$.

\end{remark}
\begin{table}
    \centering
    \renewcommand{\arraystretch}{1.5}
    \begin{tabular}{c l c}
        $K$ & $\Phi_{K,0}$ & Function\\
        \hline
        $5_2, 7_{3,5}, 9_{3,6,9,16}, 10_{128,134,142}$ & $\sum_{n=0}^\infty (-1)^n q^{\frac{n(n+1)}{2}}$ & $\psi_2(q)$ \\
         $7_2, 9_4, 9_7$ &  $\sum_{n=0}^\infty (-1)^n q^\frac{n(n+1)}{2}\sum_{m=0}^n q^{-m(n-m)} \qbinom{n}{m}_{q^{-1}}$ & $\psi_3(q)$\\
         $8_{15}, 9_{49}, 10_{66}$ & $\sum_{n=0}^\infty (-1)^n q^{\frac{n(n+1)}{2}}\sum_{m=0}^n q^{-m(n-m)} \qbinom{n}{m}_q$  & -- \\
         $12n_{293,321,750,830}$& $\sum_{n=0}^\infty q^{n^2}$ & -- \\ 
    \end{tabular}
    \caption{Closed formulas for tails of various strongly quasinegative knots.}
    \label{tab:tails-closed-formulas}
\end{table}

\subsubsection{Double twist knots} 

\begin{figure}
  \centering
  \includegraphics[scale=0.8]{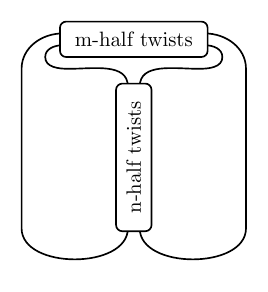}
  \caption{The double twist links $C_{m,n}$.}
\label{fig:twist_knots}
\end{figure}

We denote by $C_{m,n}$ the double twist link with $m$ and $n$ negative half twists, as depicted in \cref{fig:twist_knots}. They are two-bridge and, in particular, alternating links. The link $C_{m,n}$ is a knot if and only if $m\cdot n$ is even.

For $m,n>0$, the depicted diagrams of double twist knots $C_{2m,2n}$ and $C_{2m+1,-2n}$ have only negative crossings, so those knots are are negative, and hence strongly quasinegative \cite{Rudolph1999}. For example, $C_{2,2}$ is the negative (right-handed) trefoil.
Park proposed an explicit formula for the $F_K$ invariant for these families, based on formulas for their colored Jones polynomials \cite[Eqs. 34--39]{Park20}. If we use these formulas, we can easily extract the minimal powers.

\begin{proposition}\label{prop:L-twist-knots} For $m,n>0$,
\begin{align*}
    \el(C_{2m,2n}) &=-1 =-g(C_{2m,2n}), \\
    \el(C_{2m+1,-2n}) &=-n =-g(C_{2m+1,-2n}).
\end{align*}
\end{proposition}

It is natural to ask whether the same $F_K$ series can be obtained using the stratified sum. For the simplest non-fibered example, the negative $5_2$ knot $C_{4,2}$, the stratified sum converges and gives the same answer.

\begin{example}\label{ex:7_2}
For $C_{5,2}$, or $7_2$ in Rolfsen's table, the situation is more complicated. Park's explicit formula of $F_K$ for double twist knots gives
\begin{align*}
    f_0(q) &= -q^{-1} + 1 - q^{4} + q^{7} - q^{15} + q^{20} - q^{32} + q^{39} - q^{55} + q^{64} - q^{84} - \cdots\\
    f_1(q) &= -q^{-1} + 1 + q^3 - q^4 - q^5 - q^6 + q^7 + q^8 + q^{13} + q^{14} -q^{15} - q^{16}+\cdots
\end{align*}
These series coincide with the beginning of the stability series, as expected. If we now consider the stratified state sum, computed using the mirror of its strongly quasipositive braid $\sigma_1^2\sigma_{2,3}\sigma_{1,2}\sigma_3$, as tabulated in KnotInfo, the first coefficient converges to the series $f_0$ above. However, the second coefficient oscillates between very different expressions for odd and even weights. Those eventually stabilize to the following $q$-series:
\begin{align*}
f_1^{\text{even}}(q) &= q^{-1} + 7 + 12 q + 24 q^2 + 43 q^3 + 69 q^4 + 113 q^5 + 177 q^6 +  271 q^7  + \cdots \\
f_1^{\text{odd}}(q) &= -3q^{-1} -5 - 12 q - 24 q^2 - 41 q^3 - 71 q^4 - 115 q^5 - 179 q^6 -  269 q^7  + \cdots 
\end{align*}
We find that the \emph{average} of the above series gives the expected series $f_1$:
\[
f_1^{\text{even}}(q) + f_1^{\text{odd}}(q) = 2 f_1(q).
\] 
By analyzing the states contributing to $f_1^{\text{even}}(q)$ and $f_1^{\text{odd}}(q)$, it is possible to find a different order of summation with which the sum converges to $f_1(q)$. However, stabilizing the next term $f_2(q)$ already becomes somewhat complicated, and will require further study.

If one uses the braid representative $\sigma_1\sigma_2\sigma_4\sigma_3\sigma_4\sigma_{1,3}$ of $7_2$ instead, the state sum seems to converge for every term and matches Park's explicit formula.
\end{example}

\begin{remark} The behavior of the $F_K$ invariant under mirroring is not known. As explained above, for $C_{4,2}=5_2$, we have $\el(5_2)=-g(5_2)=-1$. In \cite{ParkThesis}, Park proposed two candidates for the $F_K$ series of the mirror knot $C_{-4,-2}$, one based on the $C$-polynomial recurrence, and another based on a certain regularization of the state sum. The first coefficient $f_0$ for both proposals are
\[
    -q^2 + q^4 - q^5 +2q^6 - 2q^7 + q^8 - 3q^{10} + \cdots,
\]
and
\[
    q^2 + 3q^3 + 6q^4 + 13 q^5 + 23q^6 + 44q^7 + 74 q^8 + \cdots,
\]
respectively. Note that, if either of the two were the correct $F_K$, we would have $\el(\mir{5_2})=2$, breaking the symmetry
$\el(\mir{K})=-\el(K)$ that holds for nice knots (\cref{prop:L-mirror}). It would also violate the formula \eqref{eq:l_hfk_intro} relating $\el$ to a grading appearing in knot Floer homology.
\end{remark}
\subsection{MMR expansion and radial limits}\label{subsec:radial_limits}

We checked that the expansion at $q=1$ of some strongly quasinegative knots recovers the MMR expansion. In particular, for all knots of at most 10 crossings appearing in \cref{tab:tails}, we checked numerically that
\begin{equation}
   \lim_{q\nearrow 1} f_0(q)\simeq-\frac{1}{\lc \Delta_K}.
\end{equation}

\begin{example} The Melvin--Morton--Rozansky expansion (in $y=x+x^{-1}-2$) of $5_2$ is
\begin{equation}\label{eq:5_2-MMR}
    \sum_{j=0}^\infty \frac{P^{(j)}_{5_2}(y)}{\Delta_{5_2}^{2j+1}(y)}h^j
\end{equation}
with $\Delta_{5_2}(y) = 1+ 2 y$ and
\begin{equation}\label{eq:5_2-MMR-polys}
\begin{split}
    P^{(1)}(y)&= 6y+5y^2, \quad P^{(2)}(y) = 2-7 y+36 y^2+54 y^3+23 y^4,\\
    P^{(3)}(y)&= 4 - 83 y + 140 y^2 - 156 y^3 - 467 y^4 - 358 y^5 - 103 y^6.
\end{split}
\end{equation}
Expanding \eqref{eq:5_2-MMR} at $x=0$ gives
\begin{equation*}
\begin{split}
    &\left(\frac{1}{2}-\frac{5 h}{2^3}+\frac{23 h^2}{2^5}-\frac{103 h^3}{2^7}
    +\dots\right)x
    +\left(\frac{3}{2^2}-\frac{17 h}{2^4}+\frac{85 h^2}{2^6}-\frac{407 h^3}{2^8}
    +\dots\right)x^2\\&+\left(\frac{5}{2^3}-\frac{15 h}{2^4}+\frac{169 h^2}{2^7}-\frac{117 h^3}{2^6}
    +\cdots\right)x^3+O(x^4)
\end{split}
\end{equation*}
which numerically agrees with 
\begin{equation*}
    \sum_{j,k=0}^\infty \left(\lim_{q\nearrow 1}\frac{d^k\tilde{f}_j(q)}{dq^k}\right) \frac{h^k}{k!} x^{j+\frac12}.
\end{equation*}
In particular, the coefficient $\tilde{f}_0(q)$ is
\begin{equation*}
    \tilde{f}_0(q)=-f_0(q)=q^{-1}\psi_{2}(q) = q^{-1}\sum_{n=0}^\infty (-1)^n q^\frac{n(n+1)}{2}, 
\end{equation*}
and one may easily check that
\begin{equation*}
    \lim_{q\nearrow 1} \tilde{f}_0(q)= \frac12,\quad \lim_{q\nearrow 1} \frac{d\tilde{f}_0(q)}{dq}=-\frac58, \quad \lim_{q\nearrow 1} \frac12\frac{d^2\tilde{f}_0(q)}{dq^2}=\frac{23}{32}.
\end{equation*}
\end{example}

\begin{example} While the $F_K$ invariant recovers the MMR expansion in many examples, this does not determine the $F_K$ coefficients. Consider the strongly quasinegative knots $12n_{148}$ and $13n_{1533}$. Although both share the same Alexander polynomial $x^{-3} + x^{-2} - 5x^{-1} + 7 - 5 x + x^2 + x^3$, the former is fibered, while the latter is not. 

We computed the $F_K$ invariant of both knots using their strongly quasinegative braid representatives, as tabulated in \cite{knotinfo} and \cite{StoimenowKnotData} respectively. The invariant associated to $12n_{148}$ has Laurent polynomials in $q$ as $x$-coefficients:
\begin{equation*}
    F_{12n_{148}}(x,q)= x^\frac12\left(-q^{-3}x^2+2q^{-3}x^3+\left(-q^{-5}-4q^{-4}-2q^{-2}\right)x^4+ O(x^5)\right);
\end{equation*}
while the $F_K$ invariant associated to $13n_{1533}$ has formal Laurent $q$-series as $x$-coefficients:
\begin{equation*}
\begin{split}
    F_{13n_{1533}}(x,q) =& x^\frac12\left(-q^{-3}\left(1-q^2+2q^4-q^6-q^7+2q^9-q^{10}-q^{12}+\cdots\right) x^2 \right. \\
    & \left. -q^{-2}\left(1 - 3 q - 2 q^2 + 4 q^3 + 2 q^4 - 2 q^5 - 3 q^6 + 5 q^8 +\cdots \right)x^3+O(x^4)\right).
\end{split}
\end{equation*}
\end{example}

We view the previous example as evidence that the finiteness of the $F_K$ coefficients could be \emph{equivalent} to the fiberedness of the knot $K$ (Park conjectured one direction in \cite[Conj.\! 1]{Park21}). This would mean that the $F_K$ invariant detects fiberedness, unlike the colored Jones function.

\subsubsection{Divergent examples}

The relation of $F_K$ to the MMR expansion does not persist for all strongly quasinegative knots. Indeed, \cref{prop:stratified-ground-state} shows that the ground state in the stratified state sum detects the genus of the knot, whereas the Alexander polynomial appearing in the MMR expansion does not necessarily do so. The strongly quasinegative knots up to 12 crossings for which this happens are $12n_{293}, 12n_{321}, 12n_{750}$ and $12n_{830}$. For all these knots, the leading term of the stratified sum is
\begin{equation*}
    -(q^{-1} x)^{\genusK}\sum_{n\geq 0} q^{n^2}= -(q^{-1}x)^\genusK (1+q+q^4+q^9+q^{16}+\cdots).
\end{equation*}
At $q \nearrow 1$, the series grows asymptotically as
\begin{equation*}
    -\frac{1}{2}q^{-\genusK}\left(1+\frac{\sqrt{\pi}}{\sqrt{-\log{q}}}\right).
\end{equation*}
These knots also satisfy that $P^{(1)}_K$ provides a better bound for the genus than the Alexander polynomial\footnote{We checked this fact using the code in \cite{BNVdV18}.}. Therefore, the MMR expansion at $x=0$ is not in $\Q\pr{x,h}$---as would be expected if we take radial limits of $F_K$---but in $\Q\lr{x}\pr{h}$.

Not surprisingly, the $F_K$ of these knots also do not recover the ADO polynomials at roots of unity. For this reason, we do not expect \cref{conj:strong-FK-ADO} to hold in general.

\begin{example} The MMR expansion of $K=12n_{293}$ in the variable $y$ is
\begin{equation*}
    \frac{1}{1+2 y}+h\frac{\left(-2 y+5 y^2+12 y^3+2 y^4\right)}{(1+2 y)^3}+O(h^2).
\end{equation*}
If we expand it in the variable $x$ at $x=0$, we get
\begin{equation*}
    \frac{x}{2}+\frac{3 x^2}{4}+\frac{5 x^3}{8}+h \left(\frac{1}{4} x^{-1}+\frac{5}{8}-x-\frac{7 x^2}{4}-\frac{69 x^3}{64}\right)+O(h^2).
\end{equation*}
The appearance of negative exponents of $x$ is related to the fact that $P^{(1)}_K$ detects the genus of the knot, unlike the Alexander polynomial:
\begin{equation*}
    \deg \Delta_{K}(x)=-1> -2 = -g(K),\quad\deg P^{(1)}_{K}(x)=-4=-2g(K).
\end{equation*}
\end{example}

\begin{example}
Another example for which the relation of $F_K$ to the MMR expansion fails is the knot $12n_{801}$. The leading term of the stratified sum reads:
\begin{equation*}
    f_0(q) = -q^{-2}(1 + 2 q + q^2 + q^3 + 2 q^4 + 2 q^5 + 2 q^6 + 2 q^7 + q^8 + 3 q^9 +\cdots).
\end{equation*}
Approaching $q=1$, $f_0$ appears to grow roughly as\footnote{We are grateful to M. Storzer for providing an improved approximation of the constant in the exponential term, compared with the first version of this paper.}
\begin{equation*}
     f_0(q) \sim -\exp{\left(-\frac{\pi^2}{50\log{q}}\right)}\frac{\sqrt{\pi}}{\sqrt{-\log{q}}}.
\end{equation*}
The Alexander polynomial of $12n_{801}$ is monic and its degree coincides (up to sign) with the genus. However, it still cannot be the Alexander polynomial of a fibered strongly quasipositive knot, since its leading term is negative.
\end{example}

\subsection{Similarities with knot Floer homology}\label{sec:connection-HFK}

The knot Floer homology $\widehat{\mathit{HFK}}(K)$ of a knot $K$,  defined by Ozsv\'ath and Szab\'o \cite{Ozsvth2004}, is a doubly-graded homology theory which categorifies the Alexander polynomial. It decomposes as
 \[
 \widehat{\mathit{HFK}}(K) = \bigoplus_{m,n \in \Z} \widehat{\mathit{HFK}}_m(K,n),
 \]
where $n$ denotes the \emph{Alexander grading} and $m$ is the \emph{Maslov (homological) grading}.

For fibered knots, $\widehat{\mathit{HFK}}(K,g)$ is one-dimensional, and the unique Maslov grading is given by the Hopf invariant. On the other hand, we saw in \cref{subsec:hopf} that the Hopf invariant appears as the minimal $q$-power $\el(K)$ of the leading coefficient of $F_K$ (\cref{thm:L-Hopf-homogeneous}, \cref{conj:L-hopf}). In the previous section, we also found a definition of $\el(K)$ for non-fibered knots. This raises the question of whether there is a global relation between $F_K$ and $\widehat{\mathit{HFK}}(K)$ that extends beyond fibered knots. Indeed, we can reformulate \cref{conj:L-hopf} as follows:
\begin{equation}\label{eq:L-HFK}
    \el(K)=g(K)-\min\{m\mid\widehat{\mathit{HFK}}_m(K,g) \neq 0\}.
\end{equation}

For strongly quasinegative knots, we have $\ell(K) = -g(K)$. Since the knot Floer homology of fibered and alternating knots is supported in a single Maslov grading, the choice of selecting the minimal Maslov grading in \eqref{eq:L-HFK} comes from the following observation about the knot Floer homology of strongly quasipositive knots:
\begin{conjecture}\label{conj:HFK-SQP} Let $K$ be a strongly quasipositive knot with genus $g$. Then the minimal Maslov grading of the top Alexander piece of $\widehat{\mathit{HFK}}$ is zero:
\begin{equation}
    \min\{m\mid\widehat{\mathit{HFK}}_m(K,g) \neq 0\} = 0. 
\end{equation}
\end{conjecture}

This holds for all non-alternating strongly quasipositive knots up to 16 crossings. If, in addition, $K$ is fibered or positive, the top Alexander piece is supported on a unique Maslov grading and the conjecture follows from \cite{HFKfibered02,Ni2006}.

The connection of the $F_K$ invariant and knot Floer homology seems to go further than just determining the invariant $\el(K)$. In \cref{subsec:radial_limits}, we mention five knots for which the stratified sum diverges at $q=1$ and does not recover the MMR expansion. Four of these knots, $12n_{293,321,750,830}$, were identified using $\dK<\genusK$. These four knots have the same leading term in $F_K$ (up to an overall $x$ and $q$ power), and their knot Floer homology satisfies
\begin{equation*}
    \widehat{\mathit{HFK}}_m(K;g) \cong\begin{cases}
        \Z \quad &\text{if } m=0 \text{ or } $1$, \\
        0 \quad &\text{otherwise.}
    \end{cases}
\end{equation*}
For $12n_{801}$, we have
\begin{equation*}
    \widehat{\mathit{HFK}}_m(K;g) \cong\begin{cases}
        \Z \quad &\text{if } m=0, \\
        \Z^2 \quad &\text{if } m=1, \\
        0 \quad &\text{otherwise.}
    \end{cases}
\end{equation*}
These five knots are the only strongly quasinegative knots up to 12 crossings (that we have identified) for which both the $F_K$ invariant diverges at $q=1$, and their knot Floer homology is not supported in a single Maslov grading. Moreover, even though $12n_{801}$ has monic Alexander polynomial, the $F_K$ series behaves very differently than that of a fibered knot. It seems that $F_K$, like knot Floer homology, is reading more refined information than the Alexander polynomial.

\begin{remark}
    We also observed that, for strongly quasipositive knots up to 16 crossings,
\begin{equation}
    \widehat{\mathit{HFK}}_m(K,n) \neq 0  \implies  n-g(K)+2 \geq m \geq n-g(K).
\end{equation}
While we do not see any reason for the upper bound to hold in general, the lower bound could somehow reflect the positivity condition.
\end{remark}

\begin{remark}\label{rmk:Murasugi}
    We believe that a possible reason behind the similarity between the extremal knot Floer homology $\widehat{\mathit{HFK}}$ and the minimal coefficient of $F_K$ could be the following.
    
    Cheng, Hedden and Sarkar in \cite{Cheng2024} defined the \emph{Grothendieck group of links under Murasugi sum along minimal index surfaces} $K(\text{links},*)$, and showed that the extremal knot homology defines a homomorphism from $K(\text{links},*)$ to the multiplicative group $\Q^{\times}_{>0}(t)$ of the ring of polynomials with positive rational coefficients. \cref{conj:HFK-SQP_intro} restricts the image of the subgroup of $K(\text{links},*)$ generated by strongly quasipositive links.
    
    We believe that the same may be true for the first coefficient $f_0(q)$ of $F_K$ (once suitably defined), namely that $\tilde{f}_0(q) = - f_0(q)$ would be multiplicative under Murasugi sum along minimal index surfaces, and consequently define a homomorphism to the ring of Laurent power series $\Z\lr{q}$. \cref{tab:tails} would then demonstrate that there are classes in $K(\text{links},*)$ that are not distinguished by extremal knot homology, but are by $f_0(q)$, answering negatively a question posed in \cite{Cheng2024}. An example would be the classes of $9_4$ and $8_{15}$; see \cref{tab:tails}.
    
    For the subgroup of $K(\text{links},*)$ generated by strongly quasipositive links\footnote{Strongly quasinegative links in our conventions.}, this homomorphism would be realized by the tail of the colored Jones polynomial. We checked this on several examples; e.g. the knot 
    $16_{592787}$ in Stoimenov's table, represented by the braid word $\sigma_4 \sigma_{1,2}\sigma_4 \sigma_3\sigma_5\sigma_{2,3}\sigma_{4,5}\sigma_{1,2}\sigma_3\sigma_2^2$, is a Murasugi sum of $5_2$ and $12n_{293}$, represented by the braid words $\sigma_1^2\sigma_2\sigma_{1,2}$ and $\sigma_{1,2}\sigma_3\sigma_{2,3}\sigma_{1,2}\sigma_3\sigma_2^2$ respectively (see \cite{Rudolph1998}). Using that
    \begin{align*}
        f_0^{5_2}(q) = -q^{-1} \sum_{n \geq 0} (-1)^n q^{\frac{n(n+1)}{2}} &= -q^{-1} (1-q+q^3+q^6 - \cdots),\\
        f_0^{12n_{293}}(q) = -q^{-2} \sum_{n \geq 0} q^{n^2} &= -q^{-2} (1+q+q^4+q^9+\cdots),    
    \end{align*}
    we checked that
    \[
    f_0^{16_{592787}}(q) = -q^{-3}(1 - q^2 + q^3 + 2 q^4 - q^5 - q^6 + q^9 + \cdots) \equiv -f_0^{5_2}(q) f_0^{12n_{293}}(q) \pmod{q^{40}}.
    \]
\end{remark}

\section{Further questions and future directions}

We end with some additional questions and directions for future work.

\medskip
\textbf{Questions:}
\begin{enumerate}
    \item Does every fibered knot admit a braid with inversion datum and is there an algorithm to find it? Is the inversion datum always unique, up to obvious symmetries? Does the inversion datum have a geometric interpretation?
    \item The Hopf invariant can be computed combinatorially both from the inversion datum on a braid (see \eqref{eq:L-general-formula}) and from a grid diagram using grid homology. Is there any relation between both models?
    \item What is special about $F_K$ when $K$ is alternating? Does it satisfy $\deg_q(f_0)=\deg_q(f_1)$? If $K$ is also fibered, is the coefficient $f_1(q)$ a monomial (this holds up to 12 crossings)?
    \item Is there a relation between the $F_K$ slopes and the exceptional slopes?
\end{enumerate}

\medskip
\textbf{Future directions:}
\begin{enumerate}
    \item Prove \cref{conj:FK-ADO} for nice fibered knots, perhaps by similar methods as \cref{thm:FK-nice-knot-MMR}.
    \item Prove the formula for the Hopf invariant in terms of colored Jones polynomials in \cref{conj:hopf-Jones_intro}.
    \item Find a general way to define a convergent state sum for strongly quasinegative knots. Show that the leading coefficient is multiplicative under Murasugi sums. Alternatively, show this for the tail of the colored Jones polynomial, without referring to $F_K$.
    \item Combine a suitably stratified sum with the inverted state sum to define $F_K$ for more general classes of knots.
    \item Extend the results of this paper to links and to $F_K^{\sl_N}$ \cite{Gruen22,GruSan2025}.
\end{enumerate}

\crefalias{section}{appendix}

\appendix

\section{Mirrors and connected sums}

In this appendix, we show that the class of nice knots is closed under taking mirror images and connected sums, and deduce the behaviour of the $F_K$ invariant under those operations.  Although this may be known to experts, we included it for a lack of reference. 

\begin{proposition}\label{prop:L-mirror} Let $K$ be a nice knot. Then, its mirror $\mir{K}$ is also a nice knot and satisfies 
$F_{\mir{K}}(x,q)=F_K(x,q^{-1})$. 
In particular, $\el(\mir{K})=-\el(K)$.
\end{proposition}

\begin{proof} Let $\bd$ be a braid diagram of $K$ that admits an inversion datum $\id$. Let $\mir{\bd}$ be the mirror image of $\bd$ (all the signs of every crossing $c$ in $\bd$ are flipped). We can naturally identify $E_\bd$ and $E_{\mir{\bd}}$ and set $\mir{\id}\coloneqq-\id$ (all signs are flipped w.r.t $\id$). Then, $\mir{\id}$ is an inversion datum for $\mir{\bd}$, since flipping the signs ensures that \eqref{eq:allowed-signs-crossings} is satisfied.

Let us first show that $(\mir{K},\mir{\bd},\mir{\id})$ is a nice knot, as in \cref{def:nice-knot}. There is a one-to-one correspondence between the set of valid states $\Omega(\id)$ and $\Omega(\mir{\id})$, given by $\mir{s}=-s-1$, with $s\in\Omega(\id)$ and $\mir{s}\in\Omega(\mir{\id})$.
In particular, we have that, for $s\in\Omega(\id)$ and the corresponding $\mir{s}\in\Omega(\mir{\id})$,
\begin{equation*}
    \frac{\mir{s}(e_{BR})+\mir{s}(e_{TR})+1}{2} = -\frac{s(e_{BR})+s(e_{TR})+1}{2}.
\end{equation*}
Therefore, for a crossing $c\in V_\bd$ and the corresponding $\mir{c}\in V_{\mir{\bd}}$,
\begin{equation*}
\begin{split}
    \dx(R^{\sgn(\mir{c})}(\mir{s})) &= \sgn(\mir{c})\frac{\mir{s}(e_{BR})+\mir{s}(e_{TR})+1}{2} \\&
    = -\sgn(c)\left(-\frac{s(e_{BR})+s(e_{TR})+1}{2}\right) = \dx(R^{\sgn(c)}(s))
\end{split}
\end{equation*}
We conclude that $\dx(P(s))=\dx(P(\mir{s}))$ for every valid state $s\in\Omega(\id)$ and its corresponding $\mir{s}\in\Omega(\mir{\id})$. In particular, since $(K,\beta,\id)$ satisfy the conditions in \cref{def:nice-knot}, $(\mir{K},\mir{\bd},\mir{\id})$ do too, and $\mir{K}$ is a nice knot.

Having proved that $\mir{K}$ is nice, the relation $F_{\mir{K}}(x,q)=F_K(x,q^{-1})$ is now a consequence of the relation between the $F_K$ series and the MMR expansion. 
In this case, it is more convenient to work with the exponential MMR expansion, instead of the integral one used in \cref{thm:FK-nice-knot-MMR}. The Habiro series evaluated at $q=e^h$ and expanded at $h=0$ satisfies\footnote{We note that the numerators $P^{(j)}_{K}$ and $Q^{(j)}_{K}$ of the integral and exponential expansions, respectively, are related by a simple change of basis.}
\begin{equation}\label{eq:Habiro-MMR-exponential}
    J_{K}(y,e^h)=\sum_{j\geq 0} \frac{Q_{K}^{(j)}(y)}{\Delta_{K}^{2j+1}(y)}h^j
\end{equation}
where $Q_{K}^{(0)}=1$, $Q_K^{(0)}\in\Q[x+x^{-1}]$ for every knot $K$ and the right-hand side is expanded at $y=0$. 

The Habiro coefficients satisfy $\an{n}{\mir{K}}=\anq{n}{K}{q^{-1}}$, which implies that $Q_{\mir{K}}^{(2j)}(y)=Q_{K}^{(2j)}(y)$ and $Q_{\mir{K}}^{(2j+1)}(y)=-Q_{K}^{(2j+1)}(y)$ for every $j\geq 0$. By \cref{thm:FK-nice-knot-MMR}, $F_{K}(x,e^h)$ recovers the expansion of the right-hand side of \eqref{eq:Habiro-MMR-exponential} at $x=0$ and $h=0$ for any nice knot $K$. Moreover, since $K$ and $\mir{K}$ are nice knots, the $x$-coefficients of $F_K(x,q)$ and $F_{\mir{K}}(x,q)$ are Laurent polynomials in $q$, and hence are uniquely determined by their expansions in $h$. Therefore, we must have $F_{\mir{K}}(x,q)=F_K(x,q^{-1})$. \qedhere
\end{proof}

\begin{remark}\label{rmk:ell-mirror}
    We can use \eqref{eq:L-general-formula} to prove the statement $\el(\mir{K})=-\el(K)$ directly. Under mirroring, the signs of the $b_i$ segments are flipped and $r_q^{\sgn(\mir{c}),\mir{\id}}=-\Rq$ (see \cref{tab:xq-contributions}). It follows that $\el(\mir{K})=-\sum_{j=2}^n\frac{\mir{\id}(b_j)}{2}+\sum_{\mir{c}\in V_{\mir{b}}} r_q^{\sgn(\mir{c}),\mir{\id}} = \sum_{j=2}^n \frac{\id(b_j)}{2}-\sum_{c\in V_\bd} \Rq =-\el(K).$
\end{remark}

\begin{proposition}\label{prop:L-connected-sum} Let $K_1$ and $K_2$ be two nice knots. Then, the connected sum $K_1\# K_2$ is a nice knot and satisfies
\begin{equation*}
    F_{K_1\# K_2}(x,q)=\frac{F_{K_1}(x,q)F_{K_2}(x,q)}{x^{\frac12}-x^{-\frac12}}.
\end{equation*}
In particular, $\el(K_1\# K_2)=\el(K_1)+\el(K_2)$.
\end{proposition}
\begin{proof} Let $\bd_i$ be a braid representative of $K_i$ with an inversion datum $\id_i$ for $i=1,2$. Let $n_1$ and $n_2$ be the number of strands of $\beta_1$ and $\beta_2$, respectively. Denote by $\sh_j\colon B_n\rightarrow B_{n+j}$ the homomorphism defined via $\sh_j(\sigma_{i})=\sigma_{i+j}$. 
Define $m\coloneqq\id_2(b_{{n_1}+1})\in\{-1,1\}$ as the sign assigned to the bottom of the leftmost strand in $\sh_{n_1}\!(\beta_2)$.

A braid representative of $K_1\#K_2$ is given by $\bd_3=\bd_1\sigma_{n_1}^{m}\mathrm{sh}_{n_1}\!(\bd_2)$. Let $\id_3$ be the function on $E_{\bd_3}$ induced by $\id_1$ and $\id_2$, where we extend to the adjacent segments of $\sigma_{n_1}^m$ the signs associated to $b_{n_1}$ and $b_{n_1+1}$ in $\bd_1$ and $\sh_{n_1}\!(\bd_2)$ respectively. For example, in \cref{fig:braid-connected-sum}, $\id_3(e_{n_1})=\id_1(b_{n_1})$.

\begin{figure}
\centering
\includegraphics[]{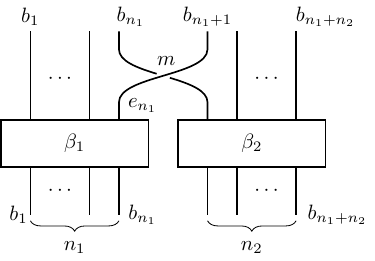}
\caption{Braid $\beta_3$ for the connected sum of the closures of $\bd_1$ and $\bd_2$.}
\label{fig:braid-connected-sum}
\end{figure}

First, let us show that $(K_1\# K_2,\bd_3,\id_3)$ satisfies \cref{def:nice-knot} and therefore is a nice knot. By definition, the signs on the segments adjacent to $\sigma_{n_1}^m$ are either $\stackanchor{++}{++}$ or $\stackanchor{-+}{-+}$ if $m=+1$, and $\stackanchor{+-}{+-}$ or $\stackanchor{--}{--}$ if $m=-1$. In particular, the minimal $x$-degree of the $R$-matrix evaluated in the segments adjacent to $\sigma_{n_1}^m$,
\begin{equation*}
    m\,\frac{s(b_{n_1})+s(b_{n_1+1})+1}{2},
\end{equation*}
is positive for every valid state $s\in\Omega(\id_3)$.

For a valid state $s\in\Omega(\id_3)$, denote by $s_i\in\Omega(\id_i)$ the associated state on $\beta_i$ for $i=1,2$. Denote by $\dx(P_1(s_1))$ the minimal $x$-degree of $(\beta_1,\id_1)$ and by $\dx(P_2(s_2))$ the minimal $x$-degree of $(\beta_2,\id_2)$ as defined in \eqref{eq:min_x_degree_P}. Then, the minimal $x$-degree of $P$ associated with $(\beta_3,\id_3)$ for an arbitrary state $s$ is
\begin{equation*}
\begin{split}
    \dx(P(s))&= \sum_{c\in V_{\beta_1}} \dx(R^{\sgn(c)}(s)) + \sum_{c\in V_{\beta_2}} \dx(R^{\sgn(c)}(s)) + m\,\frac{s(b_{n_1})+s(b_{n_1+1})+1}{2} \\&= \dx(P_1(s_1)) + \dx(P_2(s_2)) + m\,\frac{s(b_{n_1})+s(b_{n_1+1})+1}{2}.
\end{split}
\end{equation*}
Following the notation in Figure \ref{fig:braid-connected-sum}, the last equality follows from $s(e_{n_1})=s(b_{n_1})$ for any $s\in \Omega(\id_3)$ (see \cref{def:state-valid}). Therefore, since $(\beta_1,\id_1)$ and $(\beta_2,\id_2)$ satisfy the conditions in \cref{def:nice-knot}, $(\beta_3,\id_3)$ does too, and $K_1\#K_2$ is nice.

To prove the last part of the statement, the relation $\tilde{F}_{K_1\# K_2}(x,q)=\tilde{F}_{K_1}(x,q)\tilde{F}_{K_2}(x,q)$, we use the multiplicativity of the MMR expansion, which is a consequence of the multiplicativity of colored Jones polynomials $\J{n}{K_1\#K_2}=\J{n}{K_1}\J{n}{K_2}$ for every $n\geq 0$. Namely, we have
\begin{equation}
    \sum_{j\geq 0} \frac{P_{K_1\#K_2}^{(j)}(x)}{\Delta^{2j+1}_{K_1\#K_2}(x)}h^j = \left(\sum_{i\geq 0} \frac{P_{K_1}^{(i)}(x)}{\Delta^{2i+1}_{K_1}(x)}h^i\right)\left(\sum_{k\geq 0} \frac{P_{K_2}^{(k)}(x)}{\Delta_{K_2}^{2k+1}(x)}h^k\right)
\end{equation}
in the ring $\Z(x)[[h]]$. The result again follows from \cref{thm:FK-nice-knot-MMR} and the fact that the $x$-coefficients of $F_{K_1\#K_2}(x,q)$, $F_{K_1}(x,q)$ and $F_{K_2}(x,q)$ are Laurent $q$-polynomials.
\qedhere
\end{proof}

\begin{remark} As in \cref{rmk:ell-mirror}, it is easy to prove the additivity of $\ell$ using \eqref{eq:L-general-formula}. We have
\begin{equation*}
    \sum_{j=2}^{n_1+n_2} \frac{\id_3(b_j)}{2} = \sum_{j=2}^{n_1} \frac{\id_1(b_j)}{2} + \frac{\id_3(b_{n_1+1})}{2} + \sum_{j=n_1+2}^{n_2} \frac{\id_2(b_j)}{2} = \sum_{j=2}^{n_1} \frac{\id_1(b_j)}{2} + \frac{m}{2} + \sum_{j=n_1+2}^{n_2} \frac{\id_2(b_j)}{2}
\end{equation*}
and, since the $R$-matrix contribution of the crossing $\sigma_{n_1}^m$ is $\frac{m}{2}$ (\cref{tab:xq-contributions}),
\begin{equation*}
    \sum_{c\in V_{\bd_3}}\Rq = \sum_{c\in V_{\bd_1}}r_q^{\sgn(c),\id_1}+\frac{m}{2}+\sum_{c\in V_{\bd_2}}r_q^{\sgn(c),\id_2}.
\end{equation*}
Altogether, it gives $\el(K_1\# K_2)=-\sum_{j=2}^{n_1+n_2} \frac{\id_3(b_j)}{2} + \sum_{c\in V_{\bd_3}}\Rq = \el(K_1)+\el(K_2)$.
\end{remark}


\section*{References}
\printbibliography[heading=none]

\end{document}